\documentclass[fleqn,12pt]{siamart171218}
\usepackage{amsmath, amsfonts, amssymb}
\usepackage{graphicx, epstopdf}
\usepackage[export]{adjustbox}
\usepackage[caption=false]{subfig}

\begin{document}

\title{Pattern Formation in a Slowly Flattening Spherical Cap: Delayed Bifurcation \thanks{Submitted to the editors \today. \funding{This research was supported by the Natural Sciences and Engineering Research Council of Canada.}}}

\author{Laurent Charette \thanks{Institute of Applied Mathematics and Department of Mathematics, University of British Columbia, Vancouver, Canada V6T 1Z2 (\email{lchar@math.ubc.ca}).} \and Colin B. Macdonald \thanks{Institute of Applied Mathematics and Department of Mathematics, University of British Columbia, Vancouver, Canada V6T 1Z2 (\email{cbm@math.ubc.ca}).} \and Wayne Nagata \thanks{Institute of Applied Mathematics and Department of Mathematics, University of British Columbia, Vancouver, Canada V6T 1Z2 (\email{nagata@math.ubc.ca}).}}


\maketitle

\begin{abstract}
This article describes a reduction of a nonautonomous Brusselator reaction-diffusion system of partial differential equations on a spherical cap with time dependent curvature using the method of centre manifold reduction. Parameter values are chosen such that the change in curvature would cross critical values which would change the stability of the patternless solution in the constant domain case. The evolving domain functions and quasi-patternless solutions are derived as well as a method to obtain this nonautonomous normal form. The coefficients of such a normal form are computed and the reduction solutions are compared to numerical solutions.
\end{abstract}

\section{Introduction}

Pattern formation occurs in a wide array of natural phenomena. Some notable examples include coloration of animals, swarming phenomena and vortex arrangements in fluid dynamics \cite{Mein2003, Budrene1995, Chen2013}. There are also many potential applications in very active research fields such as tumour research \cite{Eftimie2017} and nanoparticle assembly \cite{Stannard2008}. Plants also exhibit different patterns through their branches, roots, leaves or petals \cite{Fujita2011, Bilsborough2011}. These patterning events are often modelled by a system of partial differential equations and the emergence of patterns can correspond to bifurcations of the solutions of the equations, where, as parameters are varied, a patternless solution loses stability and another patterned, stable solution appears. The different types of patterns can sometimes be controlled by a few parameters (size, curvature of domain, chemical concentrations or temperature for example). One specific example of pattern formation is the emergence of cotyledons in a conifer embryo. It has been observed that in the early stages of growth the embryo tip starts to flatten before cotyledons start to appear \cite{vonAd2004}.

One theory of pattern formation in morphogenesis, the process where embryonic organisms evolve into their mature shape, was formulated by Turing \cite{Turing1952} and postulates that there are chemical agents called morphogens that react among each other in a reaction-diffusion process in plant or animal tissue. The tissue then becomes differentiated and growth is then encouraged or suppressed according to the chemical pattern template. The plant hormone auxin behaves in some respects like a morphogen and is an active agent in plant growth \cite{vanMourik2012}. 

The goal of this paper is not to directly address the existence or identification of morphogens, but rather to provide a mathematical model of the transition between the patternless state and a patterned state in geometries similar to the conifer embryo. The model obtains patterns that can be similar to those observed in experiments and are often  to harmonic functions. We use a reaction-diffusion system of partial differential equations in a spherical cap domain to model pattern formation dynamics and tip geometry. More specifically we use Brusselator reaction kinetics, introduced by Prigogine and Lefever \cite{Prig1968}. This choice of kinetics in a reaction-diffusion system of equations has been studied extensively and yields bifurcation solutions that can usually be well predicted by linear analysis. It has been used previously in other pattern formation studies as well as in plant growth models \cite{Holl2008, Holl2018}.

Previous research on bifurcations of the Brusselator system in a constant spherical cap by Nagata et al. \cite{Nag2013} described a family of marginal stability curves in parameter space. The curves are the interface between the patternless and patterned regimes and can also serve to predict the type of dotted or ringed pattern emerging from the patternless solution. The different bifurcations occurring from the patternless solution in a constant spherical cap domain at one or two marginal stability curves were also studied \cite{Char2018}. In this paper we change parameters of the model in time and study the effects of this change on solutions of the nonautonomous Brusselator system.

To emulate the behaviour of the plant tip during the appearance of the physical growth patterns, we slowly decrease the curvature of the spherical cap over time while keeping the other parameters constant. We also place our initial curvature close to a constant domain marginal stability curve as we expect the transition to a patterned state should occur close to such a curve in parameter space.

One of the two main effects of the time-dependent curvature on the solutions is a slight offset of the time when the eigenvalues of the linearization of the Brusselator change sign. This is due to an extra term in the system due to the changing surface. The other effect is a delay in the effective manifestation of the bifurcation and is due to the nonautonomous nature of our equations. 

A time-dependent surface will change the Brusselator equations themselves by adding an extra term that we need to develop. Then we proceed in a similar fashion to the constant domain case, where we rewrite the equations using deviations from an underlying solution very close to the homogeneous Brusselator equilibrium, which we call the quasi-patternless solution. Next we use a projection method to reduce the nonautonomous system of partial differential equations to a centre manifold \cite{Chicone1997, Potzsche2006} and obtain a nonautonomous normal form in order to measure any delay in the bifurcation. The results are supported by numerical simulations.

\section{Model description}

We start by giving a brief description of the plant tip model used through\-out this paper. It is almost the same one used in \cite{Nag2013}, with the only difference being a time dependent curvature parameter.

\subsection{Spherical cap domains}

We model the tip of a plant embryo by a spherical cap at the end of a cylindrical stalk. A spherical cap is a part of a sphere with a circular base. One may visualize it as the surface obtained if a sphere was sliced using a horizontal plane, intersecting the plane at the base. With a fixed radius of the base $R$ and radius of the sphere $\rho = R/\gamma$, depending on a positive curvature parameter $\gamma$, the caps can range over values $0 < \gamma \leq 1$, where $\gamma = 1$ corresponds to a hemisphere and the limit $\gamma \to 0$ corresponds to a flat disk. In rectangular coordinates $(x_1,x_2,x_3)$ we align the centre of the cap on the $x_3$ axis with its circular boundary in the $x_1x_2$-plane.

\begin{figure}[ht]
	\centering
		\includegraphics[width=0.80\textwidth]{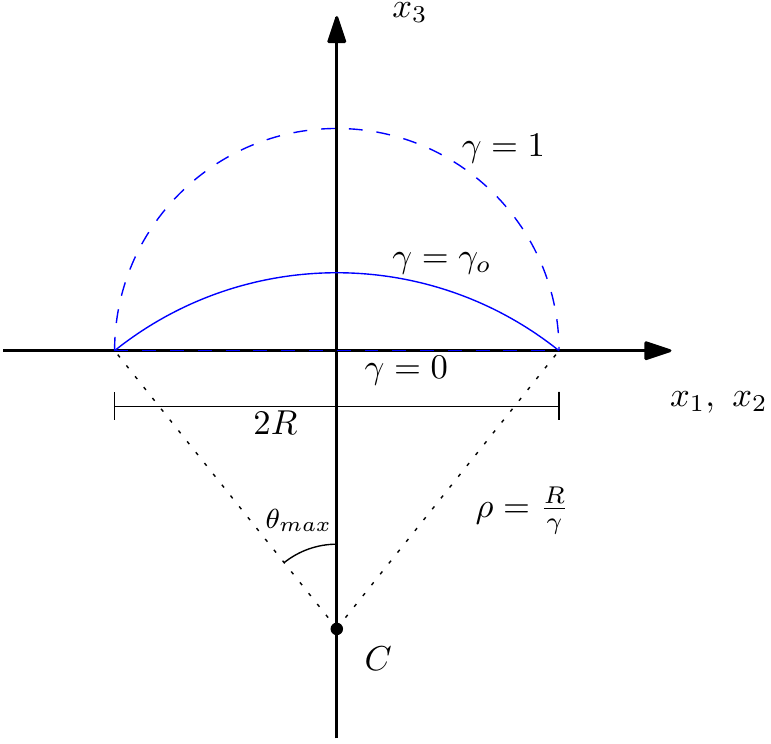}
	\caption{Cross section of a spherical cap along its diameter and the $x_3$ axis with curvature parameter $\gamma$, cap radius $R$ and spherical radius $\rho$. The spherical centre of the cap is noted by $C$. Dashed curves describe hemispherical cap ($\gamma = 0$, $C$ at the origin) and disk ($\gamma = 0$, $C$ at infinity).}
	\label{fig:sphercap}
\end{figure}

With this orientation we can use \emph{spherical coordinates} to parametrize the surface
\begin{align}
\begin{split}
x_1 = \frac{R}{\gamma} \sin \theta \cos \varphi, \ \ \  x_2 = \frac{R}{\gamma} &\sin \theta \sin \varphi, \ \ \ x_3 = \frac{R}{\gamma} \left( \cos \theta - \sqrt{1 - \gamma^2} \right),\\
\theta \in [0, \theta_{max}),& \ \ \ \varphi \in [0, 2\pi], \ \ \ \gamma \in (0, 1],
\end{split}
\end{align}
where the maximal co-latitude angle $\theta_{max}$ is dependent on the curvature $\gamma$ of the cap, and the longitudinal angle $\varphi$ is the rotation angle in the $x_1 x_2$-plane around the $x_3$-axis (also called the polar axis) computed from the $x_1$-axis. We restrict ourselves to maximal co-latitude angles up to $\pi/2$ to account for plant tip observations. In spherical coordinates the spherical cap can be described as
\begin{equation}
\Omega = \left\{ (\theta, \varphi) | \ 0 \leq \theta < \theta_{max}, \ 0 \leq \varphi \leq 2 \pi \right\}.
\end{equation}
The variable $R$ is also the radius of the cylindrical stalk and $\gamma$ is the sine of the angle the cap makes with the $x_1 x_2$-plane at the boundary. We can relate that parameter with the maximal angle using 
\begin{equation}
\theta_{max} = \arcsin \gamma.
\end{equation}
Figure \ref{fig:sphercap} illustrates how the different parameters describe the cap.

Another useful set of coordinates to describe spherical caps are \emph{toroidal coordinates}
\begin{align}\label{torcoords}
\begin{split}
x_1 = \frac{R \sinh \eta \cos \varphi}{\cosh \eta - \cos \xi}, \ \ \ x_2 &= \frac{R \sinh  \eta \sin \varphi}{\cosh \eta - \cos \xi}, \ \ \ x_3 = \frac{R \sin \xi}{\cosh \eta - \cos \xi}\\
\eta \in [0, \infty) &, \ \ \ \varphi \in [0, 2 \pi], \ \ \ \xi \in [\pi / 2, \pi).
\end{split}
\end{align}
In toroidal coordinates (\ref{torcoords}), a spherical cap is a surface with a constant $\xi$ value,
\begin{equation}
\Omega = \left\{ (\eta, \varphi) | 0 \leq \eta < \infty, 0 \leq \varphi \leq 2\pi \right\}.
\end{equation}
The variable $\xi$ measures the role of curvature here and is related to $\gamma$ by
\begin{align}
\xi = \pi - \arcsin \gamma.
\end{align}
Also worth noting is that the bounds of the toroidal coordinates $(\eta, \varphi)$ on the spherical cap do not depend on the value of $\xi$, unlike $\theta_{max}$ in spherical coordinates, which depends on $\gamma$.

\subsection{Reaction-diffusion and linear stability in a fixed spherical cap}
In this section we give a brief overview of codimension one bifurcations in a constant, fixed spherical cap domain. We will also introduce the structure and some notation that will be followed in the following sections.

In the constant domain case, the curvature $\gamma$ is a constant parameter and defines a unique, constant spherical cap domain $\Omega$. We may use $\gamma$ as a bifurcation parameter. We solve a reaction-diffusion equation in this spherical cap domain
\begin{align}
\frac{\partial \mathbf{X}}{\partial t} &= \mathbf{D} \Delta_{\Omega} \mathbf{X} + \mathbf{f}(\mathbf{X})\label{ConXvec}
\end{align}
where $\mathbf{X} = (X, Y)^\intercal$ are the surface concentrations of the active chemical species, $\mathbf{f}(\mathbf{X}) = (f(\mathbf{X}), g(\mathbf{X}))^\intercal$ give the chemical reaction kinetics, $\mathbf{D}$ is a diagonal matrix with positive constants $D_X$ and $D_Y$ on the diagonal, and $\Delta_{\Omega}$ is the Laplace-Beltrami operator on the surface $\Omega$. The Laplace-Beltrami operator can be expressed in spherical coordinates as
\begin{align}
\Delta_{\Omega} &= \frac{\gamma}{R^2 \sin \theta}\left[ \frac{\partial}{\partial \theta} \left( \sin \theta \cdot \frac{\partial}{\partial \theta} \right) + \frac{1}{\sin^2 \theta}\frac{\partial^2}{\partial \varphi^2} \right],
\end{align}
or in toroidal coordinates as
\begin{align}
\Delta_{\Omega} &= \frac{(\cosh \eta - \cos \xi)^2}{R^2 \sinh \eta} \left[ \frac{\partial}{\partial \eta} \left( \sinh \eta \frac{\partial}{\partial \eta} \right) + \frac{1}{\sinh \eta} \frac{\partial^2}{\partial \varphi^2} \right],
\end{align}
with periodic boundary conditions for $\varphi$ at $\varphi = 0, \ 2\pi$. The functions $f$ and $g$ give the chemical kinetics, for which we use the Brusselator:
\begin{align}
\mathbf{f}(\mathbf{X}) &= \begin{pmatrix} f(X, Y) \\ g(X, Y) \end{pmatrix} = \begin{pmatrix}aA - dX - bBX + cX^2Y \\ bBX - cX^2Y\end{pmatrix},\label{Bdynfvec}
\end{align}
where $a$, $b$, $c$, $d$ are positive reaction rate constants and $A$, $B$ are abundant chemical concentrations, also considered to be positive constants. The algebraic system of equations $f(X,Y) = g(X,Y) = 0$ in the Brusselator has a unique constant patternless solution
\begin{equation}\label{patternless}
\mathbf{X}_{00} = (X_{00}, Y_{00})^{\intercal} = \left(\frac{aA}{d}, \frac{bBd}{aAc}\right)^\intercal.
\end{equation}
We take Dirichlet boundary conditions, setting $\mathbf{X}$ on the boundary to be equal to this patternless solution
\begin{equation}
\mathbf{X} = \mathbf{X}_{00} \ \text{on}  \ \partial \Omega,
\end{equation}
to match the concentrations in the cylindrical stalk, where they are assumed to be constant.

The reaction function $\mathbf{f}(\mathbf{X})$ has the finite Taylor expansion about the constant patternless solution $\mathbf{X}_{00}$,
\begin{align}
\begin{split}
\mathbf{f}(\mathbf{X}) = \mathbf{f}(&\mathbf{X}_{00}) + \mathbf{f}_{\mathbf{X}}(\mathbf{X}_{00})(\mathbf{X} - \mathbf{X}_{00}) + \textstyle{\frac{1}{2}} \ \mathbf{f}_{\mathbf{XX}}(\mathbf{X}_{00})(\mathbf{X} - \mathbf{X}_{00}, \mathbf{X} - \mathbf{X}_{00})\\
&+ \textstyle{\frac{1}{6}} \ \mathbf{f}_{\mathbf{XXX}}(\mathbf{X}_{00})(\mathbf{X} - \mathbf{X}_{00}, \mathbf{X} - \mathbf{X}_{00}, \mathbf{X} - \mathbf{X}_{00}),
\end{split}
\end{align}
where the different derivatives of $\mathbf{f}$ represent linear, then symmetric bilinear and trilinear functions
\begin{align}
\mathbf{f}_{\mathbf{X}}(\mathbf{X}_{00}) &= \mathbf{K}_0, & \mathbf{f}_{\mathbf{XX}}(\mathbf{X}_{00}) &= \mathbf{B}_0, & \mathbf{f}_{\mathbf{XXX}}(\mathbf{X}_{00}) &= \mathbf{C}_0,
\end{align}
with
\begin{align}
\mathbf{K}_0 &= \begin{pmatrix} k_1 & k_2 \\ k_3 & k_4 \end{pmatrix} = 
\begin{pmatrix} bB - d & \frac{a^2 A^2 c}{d^2} \\ -bB & - \frac{a^2 A^2 c}{d^2} \end{pmatrix},\label{K0def}
\end{align}
and
\begin{align}
\mathbf{B}_0(\mathbf{U}_1, \mathbf{U}_2) &= \begin{pmatrix} 1 \\ -1 \end{pmatrix} \left[ \frac{bBd}{aA}U_1 U_2 + \frac{aAc}{d}(U_1 V_2 + V_1 U_2) \right],
\end{align}
\begin{align}
\mathbf{C}_0(\mathbf{U}_1, \mathbf{U}_2, \mathbf{U}_3) &= \begin{pmatrix} 1 \\ -1 \end{pmatrix} \frac{c}{3} \left[ U_1 U_2 V_3 + U_1 V_2 U_3 + V_1 U_2 U_3 \right],\label{C0def}
\end{align}
for
\begin{align}
\mathbf{U}_j = \begin{pmatrix} U_j \\ V_j \end{pmatrix}, \ \ \ j = 1, 2, 3.
\end{align}
We also know that the constant patternless solution satisfies
\begin{equation}\label{fXis0}
\mathbf{f}(\mathbf{X}_{00}) = \mathbf{0}.
\end{equation}

We use the deviations $\mathbf{U}$ from the constant patternless solution $\mathbf{X}_{00}$
\begin{align}
\mathbf{X} &= \mathbf{X}_{00} + \mathbf{U}, & \mathbf{U} &= (U, V)^\intercal
\end{align}
in order to study the stability of $\mathbf{X}_{00}$ and the bifurcating patterned solutions that emerges from the patternless solution, which now corresponds to $\mathbf{U} = \mathbf{0}$. The reaction-diffusion system (\ref{ConXvec}) becomes
\begin{align}
\frac{\partial \mathbf{U}}{\partial t} &= \mathbf{D} \Delta_{\Omega}\mathbf{U} + \mathbf{K}_0 \mathbf{U} + \mathbf{B}_0 (\mathbf{U}, \mathbf{U}) + \mathbf{C}_0 (\mathbf{U}, \mathbf{U}, \mathbf{U}), \label{ConUvec}
\end{align}
and the Dirichlet boundary conditions for $\mathbf{U}$ become homogeneous,
\begin{align}\label{ConUBC}
\mathbf{U} &= \mathbf{0} \text{ on } \partial \Omega.
\end{align}

The spherical cap domain and Brusselator system together have symmetries under rotation of any angle $\varphi_0$ around the $x_3$-axis and reflections through vertical planes containing the $x_3$-axis, generated by the two transformations
\begin{align}\label{Symm}
\varphi &\to \varphi + \varphi_0, & \varphi &\to -\varphi.
\end{align}
For any existing solution of the differential equation (\ref{ConXvec}) or (\ref{ConUvec}), rotations and reflections of that solution will also be solutions.

\subsubsection{Marginal stability curves}\label{subsec:MargStabCurv}
To determine the stability of the patternless solution, we linearize (\ref{ConUvec}) about $\mathbf{U} = \mathbf{0}$ to obtain
\begin{align}\label{linearization}
\frac{\partial \hat{\mathbf{U}}}{\partial t} &= \mathbf{A}_0 \hat{\mathbf{U}}, & \hat{\mathbf{U}} &= \mathbf{0} \text{ on } \partial \Omega,
\end{align}
where
\begin{align}
\mathbf{A}_0 &= \mathbf{D} \Delta_{\Omega} + \mathbf{K}_0.
\end{align}
We then use the linear stability ansatz
\begin{align}
\hat{\mathbf{U}} &= e^{\sigma t}\mathbf{U}_{0},
\end{align}
to get the eigenvalue problem for $\mathbf{A}_0$,
\begin{align}
\mathbf{A}_0 \mathbf{U}_0 &= \sigma \mathbf{U}_0, & \mathbf{U} = 0 \text{ on } \partial \Omega.\label{A0eig}
\end{align}
To solve for the eigenvalues $\sigma$ we separate variables and put
\begin{align}
\mathbf{U}_0 &= \mathbf{u}_{0, mn} \Phi_{mn}, & \mathbf{u}_{0, mn} \in \mathbb{R}^2
\end{align}
where, for $m = 0, \pm 1, \pm 2, \cdots$, $n = 1, 2, 3, \cdots$ the scalar function $\Phi_{mn}$ is an eigenfunction of the Laplace-Beltrami operator with homogeneous boundary conditions,
\begin{align}
\Delta_{\Omega} \Phi_{mn} &= -\mu_{mn} \Phi_{mn}, & \Phi_{mn} &= 0 \text{ on } \partial \Omega.
\end{align}
The eigenfunctions $\Phi_{mn}$ are given in spherical coordinates by
\begin{align}
\begin{split}
&\Phi_{mn}(\theta, \varphi) = e^{im\varphi}P^m_{\lambda_{mn}}(\cos \theta), \\
&0 \leq  \theta \leq \theta_{max}, \quad 0 \leq  \varphi < 2\pi, \label{Phisph}
\end{split}
\end{align}
or in toroidal coordinates by
\begin{align}
&\Phi_{mn}(\eta, \varphi) = e^{im\varphi}P^m_{\lambda_{mn}} \left(\frac{1 - \cosh \eta \cos \xi}{\cosh \eta - \cos \xi} \right),\\
&0 \leq  \eta \leq \infty, \quad 0 \leq  \varphi < 2\pi, \label{Phitor}
\end{align}
with $\xi = \pi - \arcsin{\gamma}$. The function $P_{\lambda}^m$ is the associated Legendre function of the first kind, of integer order $m$ and real degree $\lambda$. The value of the degree $\lambda = \lambda_{mn}$ is the $n$th positive root of the equation
\begin{align}
P^m_{\lambda}\left(\sqrt{1 - \gamma^2} \right) = 0.
\end{align}
The eigenvalues $-\mu_{mn}$ of the Laplace-Beltrami operator, corresponding to the eigenfunctions $\Phi_{mn}$, are given by
\begin{align}
-\mu_{mn} = -\lambda_{mn}(\lambda_{mn} + 1)\gamma^2 R^{-2}.
\end{align}

Then for each $m$, $n$ we obtain the algebraic eigenvalue problem in $\mathbb{R}^2$,
\begin{align}
\mathbf{A}_{0,mn} \mathbf{u}_{0,mn} &= \sigma \mathbf{u}_{0,mn},
\end{align}
where
\begin{align}
\mathbf{A}_{0,mn} &= \begin{pmatrix}
-D_X \mu_{mn} + k_1 & k_2\\
k_3 & -D_Y \mu_{mn} + k_4
\end{pmatrix} 
\end{align}
We can find the eigenvalues $\sigma$ as the roots of the characteristic equation,
\begin{align}
\sigma^2 + \sigma [(D_X + D_Y) \mu_{mn} - (k_1 + k_4)] + D_X D_Y \mu_{mn}^2 - (D_X k_4 + D_Y k_1)\mu_{mn} + k_1 k_4 - k_2 k_3 &= 0.\label{Conchareq}
\end{align}
We only consider parameter values where the solutions of (\ref{Conchareq}) are real and distinct, $\sigma = \sigma_{mn}^{\pm}$ with $\sigma_{mn}^+ > \sigma_{mn}^-$. If $\sigma_{mn}^+ < 0$ for all values of $m$ and $n$, then the patternless solution is asymptotically stable because all eigenvalues of (\ref{A0eig}) are negative and all solutions of the linearization (\ref
{linearization}) approach the patternless solution exponentially in time. If $\sigma^+_{m,n} > 0$ for some $m$, $n$ then the patternless solution is unstable.

\begin{figure}[ht]
	\centering
		\includegraphics[width=1.00\textwidth]{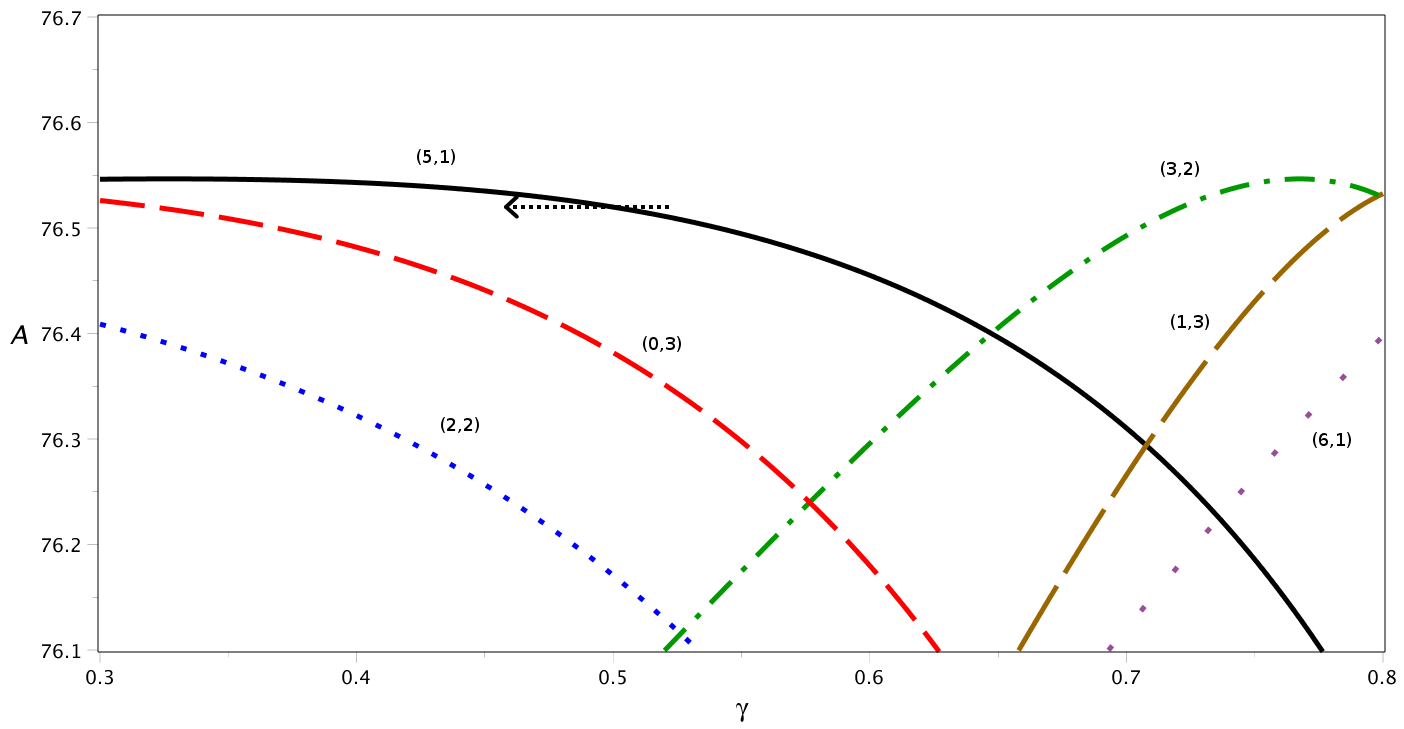}
	\caption{Marginal stability curves $A = A_{mn}(\gamma)$ for six $(m,n)$ modes using $A$ and $\gamma$ as parameters. Other parameter values are $R = 1, D_X = 0.005, D_Y = 0.1, a = 0.01, bB = 1.5, c = 1.8, d = 0.375$. These six curves are the only ones in the region with $A$ between $76.1$ and $76.7$ and $\gamma$ between $0.3$ and $0.8$, and the curves for all other $(m,n)$ values are outside the region shown. The dotted arrow shows parameter values crossing a marginal stability curve by decreasing the curvature $\gamma$ of the spherical cap through the critical value of $\gamma_0 = 0.5$, with $(m_0, n_0) = (5,1)$, $A = 76.5198$ and other parameters remaining constant.}
	\label{fig:gamvsA}
\end{figure}

By fixing all parameters except for two, we may use (\ref{Conchareq}) to obtain a marginal stability curve for each $(m, n)$ pair, where the eigenvalue $\sigma^+_{mn} = 0$. For example if we fix all parameters except for $A$ and $\gamma$ we obtain the marginal stability curve $A = A_{mn}(\gamma)$. We can compute these curves for any $(m,n)$ pair to obtain a collection of curves such as in Figure \ref{fig:gamvsA}. In this case for parameter values $(A, \gamma)$ with $A$ above all the curves, we have $\sigma^{\pm}_{mn} < 0$ for all $(m, n)$ and thus the patternless solution is stable, while if $A$ is below any one curve then $\sigma^{\pm}_{mn} > 0$ for at least one $(m, n)$ and the patternless solution is unstable. For numerical calculations we use parameter values
\begin{align}\label{param1}
R &= 1, & D_X &= 0.005, & D_Y &= 0.1,
\end{align}
\begin{align}\label{param2}
a &= 0.01, & bB &= 1.5, & c &= 1.8, & d &= 0.375.
\end{align}

Now we choose  fixed values $m_0$, $n_0$ and $\gamma_0$ so that when $\gamma = \gamma_0$ we have $\sigma^{+}_{m_0 n_0} = 0$, $\sigma^{-}_{m_0 n_0} < 0$ and $\sigma^{\pm}_{m n} < 0$ for all other $(m, n)$ pairs. Furthermore, for all $\gamma$ belonging to a sufficiently small closed interval  $[\gamma_1, \gamma_2]$ that contains $\gamma_0$ in its interior, $\sigma^{+}_{m_0 n_0}$ remains close to zero and increases from negative to positive as $\gamma$ decreases from $\gamma_2$, through $\gamma_0$, to $\gamma_1$, while all other eigenvalues remain negative and uniformly bounded away from 0, so the linearization (\ref{linearization}) has an exponential dichotomy uniformly in the parameter $\gamma$. In Figure \ref{fig:gamvsA} we show a short dotted arrow indicating parameter values that correspond to such an interval, with $\gamma_0 = 0.5$ and $A = 76.5198$ and all other parameters fixed as in (\ref{param1}) and (\ref{param2}).

\subsubsection{Bifurcation analysis}

By restricting parameter changes to within the neighbourhood of a single marginal stability curve, we can apply a codimension one bifurcation analysis. We use projection methods on the nonlinear terms to build an invariant centre manifold for the reaction-diffusion equations and obtain a normal form. This subsection briefly outlines this process for $m_0 \neq 0$.

Setting our parameters so that the largest eigenvalue of $\mathbf{A}_0$ is critical, achieved with $\sigma^+_{m_0 n_0} = 0$ for some $(m_0, n_0)$, $m_0 \neq 0$, we find its associated eigenfunction
\begin{align}
\mathbf{U}^{(0)} &= \begin{pmatrix} u_{m_0 n_0}^{(0)} \\ v_{m_0 n_0}^{(0)} \end{pmatrix} \Phi_{m_0 n_0}
\end{align}
for the critical $(m_0, n_0)$ values. This gives a two-dimensional (if $m_0 \neq 0$) centre subspace for the linearization
\begin{align}
E^c &= \left\{ z \mathbf{U}^{(0)} + \bar{z} \bar{\mathbf{U}}^{(0)} \middle| \ z \in \mathbb{C} \right\},
\end{align}
where the overbar denotes the complex conjugate and $\mathbf{U}^{(0)}$ is any fixed choice of the eigenfunction. This space denotes the dominant mode that will emerge at the bifurcation.

Because the largest eigenvalue $\sigma^{(0)} = \sigma^+_{m_0 n_0}$ is close to zero in the neighbourhood of the parameter space, there will be, for $|z|$ sufficiently small, a locally invariant centre manifold
\begin{align}
W^c &= \left\{ z \mathbf{U}^{(0)} + \bar{z} \bar{\mathbf{U}}^{(0)} + \mathcal{O}(|z|^2) \middle| \ z \in \mathbb{C} \right\}
\end{align}
for the nonlinear equation (\ref{ConUvec}), (\ref{ConUBC}). All nearby solutions of (\ref{ConUvec}), (\ref{ConUBC}) will decay exponentially to this manifold and the general long term local behaviour of the system can essentially be described by the system restricted to $W^c$. This leads us to an autonomous differential equation that describes the time evolution of the system restricted to $W^c$, a normal form for a pitchfork bifurcation
\begin{align}\label{ConNFCom}
\dot{z} &= \sigma^{(0)} z + C |z|^2 z + \mathcal{O}(|z|^5),
\end{align}
where $\sigma^{(0)}$ is the real critical eigenvalue, depending on $\gamma$ and close to $0$ that acts as the bifurcation parameter and $C$ is a coeffcient. We have that $\sigma^{(0)} = 0$ when parameter values lie on a marginal stability curve with $\gamma = \gamma_0$. This normal form must be equivariant under the transformations
\begin{align}\label{SCsym2}
z &\to e^{im \varphi_0}z, & z &\to \bar{z},
\end{align}
associated to rotation by an angle $\varphi_0$ and reflection symmetries (\ref{Symm}) respectively, and this implies that the cubic term $C$ in (\ref{ConNFCom}) is real valued.

\begin{figure}[ht]
	\centering
		\includegraphics[width=1.00\textwidth]{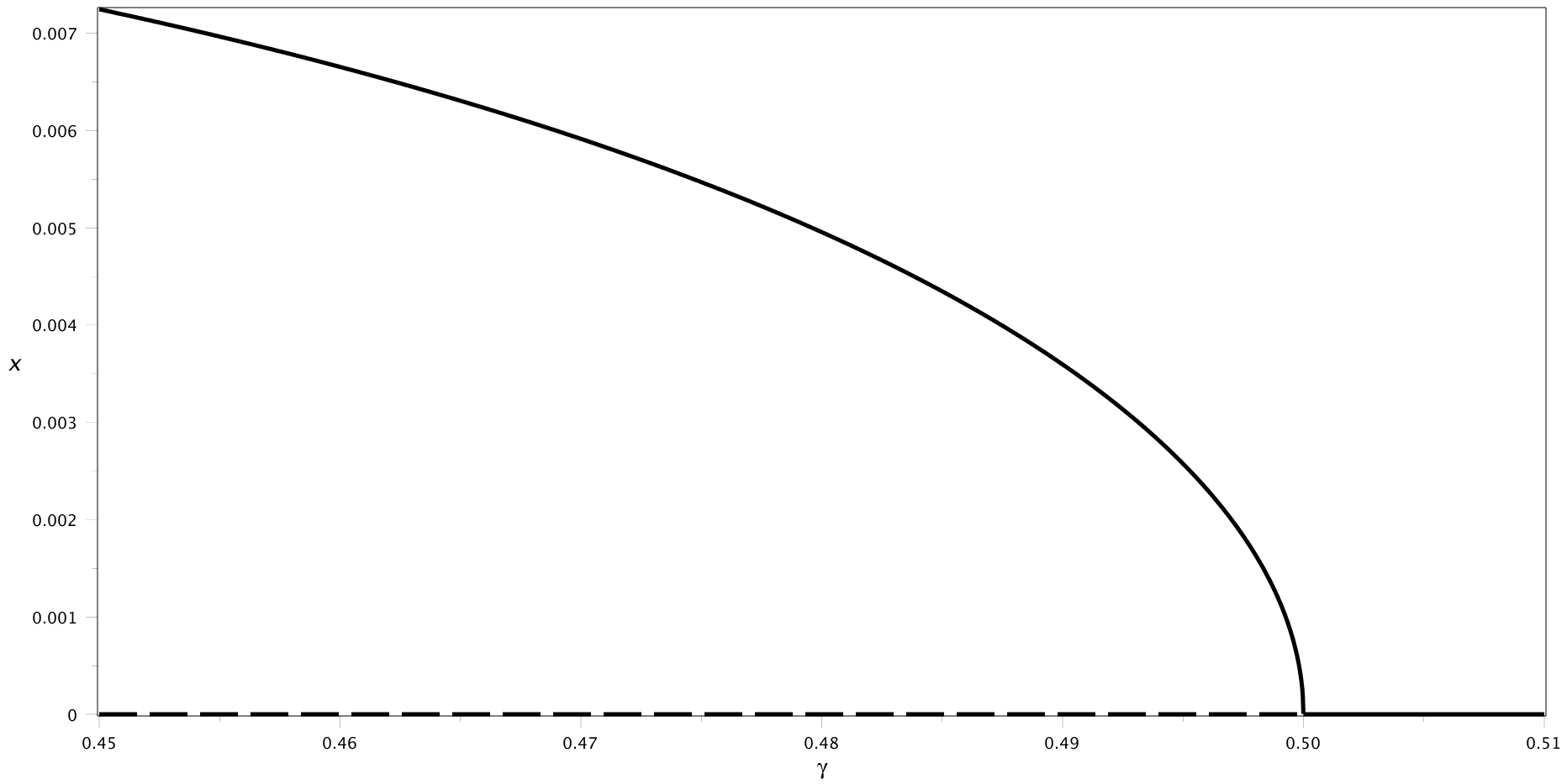}
	\caption{Bifurcation diagram of the normal form (\ref{NFxfirst}) with cubic coefficient value $C = -2.99378$. The solid curves represent stable equilibrium points and the dashed line represents the unstable equilibrium.}
	\label{fig:pitchfork}
\end{figure}

Because of this invariance under rotations, we need only consider the real part of the normal form
\begin{align}\label{NFxfirst}
\dot{x} = \sigma^{(0)} x + Cx^3 + \mathcal{O}(x^5),
\end{align}
where $x$ is the real part of $z$, and $\sigma^{(0)}$ and $C$ are the same as in (\ref{ConNFCom}). This normal form equation gives us the magnitude of the critical component of the bifurcating solution so we can locally obtain a good approximation of the corresponding solution of the partial differential equation (\ref{ConUvec}), (\ref{ConUBC}). In the case where the parameter values are given by equations (\ref{param1}, \ref{param2}) and $A = 76.5198$ the transition occurs when $\gamma = \gamma_0 = 0.5$ (as pictured by the dotted arrow in Figure \ref{fig:gamvsA}), we may solve for $\sigma^{(0)}$ as a function of $\gamma$ by taking the largest root of equation (\ref{Conchareq}). The coefficient $C$ also depends on $\gamma$, but this dependence only has a higher order effect on the normal form (\ref{NFxfirst}), and to leading order at the critical curvature $\gamma_0$ we find $C = -2.99378$. The calculation of $C$ follows the usual procedure of centre manifold and normal form reductions with symmetry (e.g. \cite{Nag2013} for more details; also cf. section 3 of this paper). The resulting bifurcation diagram is shown in Figure \ref{fig:pitchfork}.

\subsection{Reaction-diffusion in a slowly changing spherical cap}

In order to better model the flattening tip in plant embryos we make the curvature parameter $\gamma$ a function of time. This allows the spherical cap to change curvature while retaining a fixed radius $R$ at its boundary. For a slowly changing domain we will want this curvature parameter to change at a rate of order $\varepsilon$, where $\varepsilon$ is a small parameter. We can then define a slow time $\tau = \varepsilon t$ and set
\begin{equation}
\gamma = \gamma(\tau) = \gamma(\epsilon t).
\end{equation}
As a result, the time derivative becomes $\frac{\partial}{\partial t} = \varepsilon \frac{\partial}{\partial \tau}$.

As we change curvature over time the surface area of our domain also changes. This change of surface area adds extra terms in the system of equations in order to satisfy the Reynolds transport theorem \cite{Reyn1903} and maintain conservation of total chemical agents. Nonautonomous reaction-diffusion equations for chemical reactions on a growing surface $\Omega = \Omega(\tau) = \Omega(\varepsilon t)$ are
\begin{align}
\frac{\partial X}{\partial t} &= D_X \Delta_{\Omega(\varepsilon t)} X - (\nabla_{\Omega(\varepsilon t)} \cdot \mathbf{v})X + f(X, Y)\label{GroX}\\
\frac{\partial Y}{\partial t} &= D_Y \Delta_{\Omega(\varepsilon t)} Y - (\nabla_{\Omega(\varepsilon t)} \cdot \mathbf{v})Y + g(X, Y),
\label{GroY}\end{align}
where $f$ and $g$ are the Brusselator kinetics (\ref{Bdynfvec}), $\mathbf{v}$ is the velocity vector of a point on our surface, $\nabla_{\Omega(\varepsilon t)}$ is the tangential gradient to the surface. In the constant domain, where $\gamma$ is constant, we have $\mathbf{v} = \mathbf{0}$ in (\ref{GroX}), (\ref{GroY}) and in this case the solution (\ref{patternless}) represents a homogeneous equilibrium solution of the reaction-diffusion system. The Laplace-Beltrami operator $\Delta_{\Omega(\tau)}$ is now dependent on the slow time $\tau = \varepsilon t$,
\begin{align}
\Delta_{\Omega(\tau)} &= \frac{\gamma(\tau)}{R^2 \sin \theta}\left[ \frac{\partial}{\partial \theta} \left( \sin \theta
 \frac{\partial}{\partial \theta} \right) + \frac{1}{\sin^2 \theta}\frac{\partial^2}{\partial \varphi^2} \right]
\end{align}
in spherical coordinates, or
\begin{align}
\Delta_{\Omega(\tau)} &= \frac{(\cosh \eta - \cos \xi(\tau))^2}{R^2 \sinh \eta} \left[ \frac{\partial}{\partial \eta} \left( \sinh \eta \frac{\partial}{\partial \eta} \right) + \frac{1}{\sinh \eta} \frac{\partial^2}{\partial \varphi^2} \right]
\end{align}
in toroidal coordinates. We use the same equilibrium boundary values as in the constant domain case
\begin{equation}\label{XYBC}
(X, Y) = (X_{00}, Y_{00}) \ \text{on}  \ \partial \Omega(\tau).
\end{equation}
The reasoning for this boundary condition is that the base of the spherical cap remains constant over time and thus should have the same value as in the constant case.

The form of the divergence term will depend on how exactly we want the surface to change curvature. Here we choose to have points on the surface travel through points on a circle where a coordinate in toroidal coordinates is constant ($\eta$ in this case). This has the advantages of being able to be expressed by a fixed domain, having the surface velocity always normal to the surface and preserving the ratio of distances between a point and the two foci, located at each end of the bisection of the spherical cap the point in question. In general, a time-dependent surface domain can be described as a continuous mapping $\psi_{t}$ of the initial surface
$$
\Omega(\tau) = \psi_{t} (\Omega(0)) = \Omega(\varepsilon t).
$$
Under these circumstances this divergence factor can be expressed with the use of scale factors, as done by Plaza et al. \cite{Plaza2004}, for example. If our surface is described by two variables, associated with the scale factors $h_1$ and $h_2$, then
$$
\nabla_{\Omega(\varepsilon t)} \cdot \mathbf{v} = \frac{\partial}{\partial t}[\ln(h_1 h_2)].
$$

Using toroidal coordinates to compute the scale factors can greatly simplify the computations. If we keep using the spherical coordinates as before, the bounds on the polar angle $\theta$ will change as we change the curvature and finding an appropriate trajectory for a point on the moving surface becomes very cumbersome. Consequently we use toroidal coordinates and fix our parametrization of the spherical cap to $\varphi$ and $\eta$, both having the same bounds for any given curvature.

If one expresses the position $\mathbf{x} = (x_1, x_2, x_3)^\intercal$ on the surface in toroidal coordinates (\ref{torcoords}), the scale factors are given by
\begin{align}
h_1 &= \left| \frac{\partial \mathbf{x}}{\partial \eta} \right| = \frac{R}{\cosh \eta - \cos \xi(\varepsilon t)}\\
h_2 &= \left| \frac{\partial \mathbf{x}}{\partial \varphi} \right| = \frac{R \sinh{\eta}}{\cosh \eta - \cos \xi(\varepsilon t)}.
\end{align}

The divergence term may thus be computed:
\begin{equation}
\nabla_{\Omega(\tau)} \cdot \mathbf{v} = \frac{\partial}{\partial t}[\ln(h_1 h_2)] = \frac{-2 \varepsilon \xi' \sin \xi}{\cosh \eta - \cos \xi} = \frac{-2 \varepsilon z \xi'}{R},
\end{equation}
where the prime represents the $\tau$ derivative. We can then use the formulation of the $x_3$ coordinate to transform back to spherical coordinates:
\begin{equation}
\nabla_{\Omega(\tau)} \cdot \mathbf{v} = \frac{2 \varepsilon \gamma'}{\gamma} \left(\frac{\cos \theta}{\sqrt{1 - \gamma^2}} - 1 \right).
\end{equation}
We note here that when using spherical coordinates the $\theta$-coordinate will have a time dependent upper bound.

We find that the divergence term is proportional to the $x_3$-coordinate, which will depend both explicitly on time and on the distance from the point on the cap to the $x_3$-axis. We thus obtain the following equations in terms of the slow time variable $\tau$:
\begin{align}
\begin{split}
\varepsilon \frac{\partial X}{\partial \tau} &= D_X \Delta_{\Omega(\tau)} X - \varepsilon \gamma'(\tau) \ Q(\tau) X + f(X, Y)\\
\varepsilon \frac{\partial Y}{\partial \tau} &= D_Y \Delta_{\Omega(\tau)} Y - \varepsilon \gamma'(\tau) \ Q(\tau) Y + g(X, Y),
\end{split}\label{brussgro}
\end{align}
where $Q(\tau)$ depends on the position on the cap and may be expressed either in spherical coordinates by 
\begin{equation}
Q(\tau)(\theta) = \frac{2}{\gamma(\tau)} \left( \frac{\cos \theta}{\sqrt{1 - \gamma(\tau)^2}} - 1 \right), \quad 0 \leqslant \theta < \theta_{max}(\tau),
\end{equation}
with $\theta_{max}(\tau) = \arcsin \gamma(\tau)$, or in toroidal coordinates by
\begin{equation}
Q(\tau)(\eta) = \frac{2 \varepsilon \tan \xi(\tau)}{\cosh \eta - \cos \xi(\tau)}, \quad 0 \leqslant \eta < \infty,
\end{equation}
where $\xi(\tau) = \pi - \arcsin \gamma(\tau)$. We observe that the function $Q(\tau)$ does not depend on the longitudinal angle $\varphi$.

We write the system of equations (\ref{brussgro}) in vector form
\begin{equation}\label{RDEvec}
\varepsilon \frac{\partial \mathbf{X}}{\partial \tau} = \mathbf{D} \Delta_{\Omega(\tau)} \mathbf{X} - \varepsilon \gamma'(\tau) Q(\tau) \mathbf{X} + \mathbf{f}(\mathbf{X}),
\end{equation}
with boundary conditions
\begin{align}\label{RDEBC}
\mathbf{X} = \mathbf{X}_{00} \text{ on } \partial \Omega(\tau), 
\end{align}

Because the nonautonomous modifications to the equations are independent of the longitudinal angle $\varphi$, the reaction-diffusion system with our choice of domain evolution preserves the rotation and reflection symmetries (\ref{Symm}) and these will feature in the derivation of the reduced normal forms.

\subsection{The quasi-patternless solution}

Under the slow domain evolution the concentrations $X$, $Y$ change slightly so that the constant patternless solution of the autonomous constant domain system is no longer a solution of the nonautonomous system (\ref{RDEvec}). This introduces the idea of a ``quasi-patternless'' solution that will depend on time. Due to the slow curvature change, we use the slow timescale $\tau$.

We assume there is a slowly evolving, radially symmetric quasi-patternless solution $\mathbf{X}_0(\tau, \varepsilon)$ that is $\mathcal{O}(\varepsilon)$-close to the constant domain patternless solution, and coincides with it when $\varepsilon = 0$. Such a solution would have an asymptotic series
\begin{equation}\label{QPepsseries}
\mathbf{X}_0(\tau, \varepsilon) = \mathbf{X}_{00} + \varepsilon \mathbf{X}_{01}(\tau) + \mathcal{O}(\varepsilon^2),
\end{equation}
where each $\mathbf{X}_{0j}$ vector has an $X$ and $Y$ component, and the leading order term $\mathbf{X}_{00}$ is the constant domain patternless solution (\ref{patternless}).




We reduce the problem to solving for the higher order correction terms in $\varepsilon$ for the solution. We start by substituting expansion (\ref{QPepsseries}) into (\ref{RDEvec}, \ref{RDEBC}) to obtain
\begin{equation}
\varepsilon^2 \frac{\partial \mathbf{X}_{01}}{\partial \tau} + \mathcal{O}(\varepsilon^3) = \varepsilon \mathbf{D} \Delta_{\Omega(\tau)} \mathbf{X}_{01} - \varepsilon \gamma'(\tau)Q(\tau) \mathbf{X}_{00} + \mathbf{f}(\mathbf{X}_{00} + \varepsilon \mathbf{X}_{01} + \mathcal{O}(\varepsilon^2)) + \mathcal{O}(\varepsilon^2).
\end{equation}
The order $\varepsilon$ terms must satisfy
\begin{align}\label{QPfirstorderDE}
0 &= \mathbf{D} \Delta_{\Omega(\tau)} \mathbf{X}_{01}(\tau) - \varepsilon \gamma'(\tau) Q(\tau) \mathbf{X}_{00} + \mathbf{K}_0 \mathbf{X}_{01}(\tau), & \mathbf{X}_{01} &= \mathbf{0} \text{ on } \partial \Omega(\tau),
\end{align}

We can then express the solution to (\ref{QPfirstorderDE}) as a series of the now $\tau$ dependent eigenfunctions of the Laplace-Beltrami operator $\Phi_{mn} (\tau)$, where
\begin{align}\label{eigftau}
\Delta_{\Omega(\tau)} \Phi_{mn}(\tau) &= -\mu_{mn}(\tau) \Phi_{mn}(\tau), & \Phi_{mn}(\tau) &= 0 \text{ on } \partial \Omega(\tau) \text{ for all } \tau.
\end{align}

We use a series expression for the first correction term
\begin{equation}\label{QPseries}
\mathbf{X}_{01} (\tau) = \sum_{n=1}^{\infty}{\begin{pmatrix}\alpha_{0n} (\tau) \\ \beta_{0n} (\tau) \end{pmatrix} \Phi_{0n} (\tau)}.
\end{equation}
Since (\ref{QPfirstorderDE}) is independent of $\varphi$, we only write the series for order $m=0$. This makes the quasi-patternless solution invariant under rotations of the angle $\varphi$ and reflections.

To find an approximation of the quasi-patternless solution we solve for the first several coefficients $\alpha_{0n}$ and $\beta_{0n}$, and again for an increasing number of $n$ to check for convergence. We found that computing modes from $n=1$ to a few after the most unstable $m=0$ mode to be sufficient. We computed coefficients in the case when when the $(0,3)$ mode is the most unstable and the difference between computing 5 modes and computing 12 modes is very small. We compare our truncated series approximation of the quasi-patternless solution with results of a numerical computation in section 3.6.

We also expand the scalar function $Q(\tau)$ in an eigenfunction series
\begin{equation}
Q(\tau) = \sum_{n=1}^{\infty} {q_{n} (\tau) \Phi_{0n} (\tau)}.
\end{equation}
We find the coefficients $q_{n}$ by applying a projection with a Legendre function in spherical or toroidal coordinates, for example in spherical coordinates
\begin{equation}\label{qnproj}
q_{n} (\tau) = \frac{\int_0^{\arcsin \gamma(\tau)}{Q(\tau)(\theta) P^0_{\lambda_{0,n}(\tau)}(\cos \theta) \sin \theta} \ d \theta}{\int_0^{\arcsin \gamma(\tau)}{\left(P^0_{\lambda_{0,n}(\tau)}(\cos \theta) \right)^2 \sin \theta \ d \theta}}.
\end{equation}
In the end, to solve for the coefficients $\alpha_{0n}$ and $\beta_{0n}$, we need to find the solution of the linear equation
\begin{equation}\label{alphalin}
\begin{pmatrix} 0\\ 0 \end{pmatrix} = -\gamma'(\tau) q_n (\tau) \begin{pmatrix}X_{00}\\ Y_{00} \end{pmatrix} + \begin{pmatrix} -D_X \mu_{0n}(\tau) + k_1 & k_2\\ k_3 & -D_Y \mu_{0n}(\tau) + k_4 \end{pmatrix} \begin{pmatrix} \alpha_{0n}(\tau)\\ \beta_{0n}(\tau) \end{pmatrix},
\end{equation}
which is 
\begin{equation}\label{alphalinsol}
\begin{pmatrix}\alpha_{0n}(\tau)\\ \beta_{0n}(\tau) \end{pmatrix} = \gamma'(\tau) q_n(\tau) \mathbf{A}^{-1}(\tau) \begin{pmatrix}X_{00}\\ Y_{00} \end{pmatrix},
\end{equation}
where $\mathbf{A}(\tau)$ is the matrix in the right side of (\ref{alphalin}). We assume here that a critical eigenvalue $0$ of $\mathbf{A}(\tau)$ does not occur for any $(0,n)$ mode for $n \geq 1$, for any $\tau$, so the equations are all solvable. If required, we could solve for correction terms of higher order in $\varepsilon$.

\subsection{Deviation equations}
In this section we re-express the slowly changing system (\ref{RDEvec}, \ref{RDEBC}) in terms of the deviations from the quasi-patternless solution using the same methodology as in the constant domain case.


We define the deviation $\mathbf{U}$ from the quasi-patternless solution, by
\begin{equation}
\mathbf{X} = \mathbf{X}_{0}(\tau, \varepsilon) + \mathbf{U},
\end{equation}
where $\mathbf{X}_0 = (X_0, Y_0)^\intercal$ is the quasi-patternless solution and $\mathbf{U} = (U, V)^\intercal$. 
Substituting the expansion (\ref{QPepsseries}) in the system (\ref{RDEvec}, \ref{RDEBC}) with the use of expressions (\ref{K0def}) to (\ref{fXis0}) and keeping terms up to first order in $\varepsilon$ we get
\begin{align}
\frac{\partial \mathbf{U}}{\partial t} &= \mathbf{A}(\varepsilon t, \varepsilon)\mathbf{U} + \mathbf{B}(\varepsilon t, \varepsilon)(\mathbf{U}, \mathbf{U}) + \mathbf{C}_0 (\mathbf{U}, \mathbf{U}, \mathbf{U}), & \mathbf{U} &= \mathbf{0} \text{ on } \partial \Omega(\varepsilon t),\label{ABCeq}
\end{align}
where $\mathbf{A}$, $\mathbf{B}$ are respectively nonautonomous linear and bilinear operators defined by
\begin{align}
\mathbf{A}(\tau, \varepsilon)\mathbf{U} = \mathbf{D} \Delta_{\Omega(\tau)} \mathbf{U} + \mathbf{K}(\tau, \varepsilon)\mathbf{U} - \varepsilon \gamma'(\tau) Q(\tau)\mathbf{U}, 
\end{align}
\begin{align}
\mathbf{B}(\tau, \varepsilon)(\mathbf{U}_1, \mathbf{U}_2) &= \mathbf{B}_0(\mathbf{U}_1, \mathbf{U}_2) + 3 \mathbf{C}_0(\mathbf{X}_0(\tau, \varepsilon) - \mathbf{X}_{00}, \mathbf{U}_1, \mathbf{U}_2),
\end{align}
with 
\begin{align}
\mathbf{X}_0(\tau, \varepsilon) = \mathbf{X}_{00} + \varepsilon \mathbf{X}_{01}(\tau) + \mathcal{O}(\varepsilon^2),
\end{align}
\begin{align}
\begin{split}
&\mathbf{K}(\tau, \varepsilon) = \mathbf{K}_0 + 2\mathbf{B}_0(\mathbf{X}_0(\tau, \varepsilon) - \mathbf{X}_{00}, \cdot) + 3 \mathbf{C}_0(\mathbf{X}_0(\tau, \varepsilon) - \mathbf{X}_{00}, \mathbf{X}_0(\tau, \varepsilon) - \mathbf{X}_{00}, \cdot)\\
& \hspace{31px} = \mathbf{K}_0 + \varepsilon \mathbf{K}_1(\tau) + O(\varepsilon^2),\\
\end{split}\\
&\mathbf{K}_1(\tau) = \begin{pmatrix} \frac{2bBd}{aA}X_{01}(\tau) + \frac{2aAc}{d}Y_{01}(\tau) & \frac{2aAc}{d} X_{01}(\tau)\\ -\frac{2bBd}{aA}X_{01}(\tau) - \frac{2aAc}{d}Y_{01}(\tau) & -\frac{2aAc}{d} X_{01}(\tau) \end{pmatrix},\\
\begin{split}
&3\mathbf{C}_0(\mathbf{X}_0(\tau, \varepsilon) - \mathbf{X}_{00}, \mathbf{U}_1, \mathbf{U}_2) = \begin{pmatrix} 1 \\ -1 \end{pmatrix} \left[ \left( \textstyle{\frac{bBd}{aA}} + \varepsilon c Y_{01}(\tau) + \mathcal{O}(\varepsilon^2) \right) U_1 U_2\right.\\
& \quad \quad \quad \quad \left. + \left( \textstyle{\frac{aAc}{d}} + \varepsilon c X_{01}(\tau) + \mathcal{O}(\varepsilon^2) \right) \right] (U_1 V_2 + U_2 V_1).
\end{split}
\end{align}

\section{Pattern formation in a slowly changing domain}

In this section we describe emerging patterned solutions using a bifurcation theory method adapted to this slowly evolving domain case. We start by analyzing the linearization of the system (\ref{ABCeq}) to determine the stability of the quasi-patternless solution. We then perform a nonautonomous centre manifold reduction of the nonlinear system with a projection method to build our nonautonomous normal form. Our analysis includes the effects of the non-isotropically evolving domain.

\subsection{Linearization about the quasi-patternless solution}

The linearization of (\ref{ABCeq}) is
\begin{align}
\frac{\partial \hat{\mathbf{U}}}{\partial t} &= \mathbf{A}(\varepsilon t, \varepsilon) \hat{\mathbf{U}}, & \hat{\mathbf{U}} &= \mathbf{0} \text{ on } \partial \Omega(\varepsilon t).\label{Lineq}
\end{align}
Following WKB theory, e.g. \cite{Erneux1991}, we assume the following ansatz
\begin{equation} \label{WKB}
\hat{\mathbf{U}}(\tau, \varepsilon) = e^{\Psi(\tau, \varepsilon)/\varepsilon} \ \mathbf{U}(\tau, \varepsilon),
\end{equation}
to obtain
\begin{align}
\varepsilon \mathbf{U}'(\tau, \varepsilon) &= \mathbf{A}(\tau, \varepsilon) \mathbf{U}(\tau, \varepsilon) - \Psi'(\tau, \varepsilon)\mathbf{U}(\tau, \varepsilon), & \mathbf{U} &= \mathbf{0} \text{ on } \partial \Omega(\tau).\label{WKBlin}
\end{align}
Expanding in asymptotic series
\begin{align}
\begin{split}\label{U0series}
\mathbf{A}(\tau, \varepsilon) &= \mathbf{A}_{0}(\tau) + \varepsilon \mathbf{A}_{1}(\tau) + \mathcal{O}(\varepsilon^2),\\
\Psi(\tau, \varepsilon) &= \Psi_0(\tau) + \varepsilon \Psi_1(\tau) + \mathcal{O}(\varepsilon^2),\\
\mathbf{U}(\tau, \varepsilon) &= \mathbf{U}_{0}(\tau) + \varepsilon \mathbf{U}_{1}(\tau) + \mathcal{O}(\varepsilon^2), 
\end{split}
\end{align}
then substituting back into (\ref{WKBlin}), we obtain that the leading order terms must solve
\begin{align}
\mathbf{A}_{0}(\tau) \mathbf{U}_{0}(\tau) &= \Psi_0'(\tau) \mathbf{U}_{0}(\tau), & \mathbf{U}_0(\tau) = \mathbf{0} \text{ on } \partial \Omega(\tau).
\end{align}
For each $\tau$, this is an eigenvalue problem for the leading order term $\mathbf{U}_0$. If we set
\begin{align}\label{U0Phi}
\mathbf{U}_{0}(\tau) = \mathbf{u}_{0, mn}(\tau) \Phi_{mn}(\tau),
\end{align}
where $\mathbf{u}_{0, mn}(\tau) \in \mathbb{R}^2$, and $\Phi_{mn}(\tau)$ are as defined in (\ref{eigftau}), we obtain an algebraic eigenvalue problem for $\mathbf{u}_{0, mn}(\tau)$
\begin{align}\label{linearizationtau}
\mathbf{A}_{0}(\tau) \mathbf{u}_{0, mn}(\tau) &= \Psi_0'(\tau) \mathbf{u}_{0, mn}(\tau),
\end{align}
where
\begin{equation}
\mathbf{A}_{0, mn}(\tau) = \begin{pmatrix} -D_X \mu_{mn}(\tau) + k_1 & k_2\\
                               k_3 & -D_Y \mu_{mn}(\tau) + k_4\end{pmatrix}.
\end{equation}
We obtain that $\Psi'_0(\tau) = \sigma_{mn}^{\pm}(\tau)$ is one of the two (now $\tau$-dependent) eigenvalues of the matrix $\mathbf{A}_{0,mn}(\tau)$, so it solves the characteristic polynomial (\ref{Conchareq}), with $\mu_{mn} = \mu_{mn}(\tau)$ depending on the slow time variable.
Then $\mathbf{u}_{0,mn}(\tau)$ is the associated eigenvector
\begin{equation}\label{u0eigv}
\mathbf{u}_{0,mn}(\tau) = \begin{pmatrix} k_2\\ D_X \mu_{mn}(\tau) - k_1 + \sigma^{\pm}_{mn}(\tau) \end{pmatrix}.
\end{equation}

The order $\varepsilon$ terms in (\ref{WKBlin}) give us
\begin{align}
 [\mathbf{A}_0(\tau) - \Psi'_0(\tau) \mathbf{I}] \mathbf{U}_{1}(\tau) &= \Psi'_1(\tau) \mathbf{U}_{0}(\tau) + \mathbf{U}'_{0}(\tau) - \mathbf{A}_1(\tau) \mathbf{U}_{0}(\tau),& \mathbf{U}_1(\tau) &= \mathbf{0} \text{ on } \partial \Omega(\tau).\label{U1Psipeq}
\end{align}
with $\Psi'_0(\tau) = \sigma^{\pm}_{mn}(\tau)$. We assume $\mathbf{U}_1$ has a series expansion
\begin{equation}
\mathbf{U}_1(\tau) = \sum_{i=1}^{\infty} \mathbf{u}_{1, mi}(\tau) \Phi_{mi}(\tau),
\end{equation}
then (\ref{U1Psipeq}) becomes
\begin{equation}\label{u1syseq}
\left[ \mathbf{A}_{0, mi}(\tau) - \Psi'_0(\tau) \mathbf{I} \right] \mathbf{u}_{1, mi}(\tau) = \delta_{ni} \left[ \Psi'_1(\tau) \mathbf{u}_0(\tau) + \mathbf{u}'_{0, mn}(\tau) \right] + d_i(\tau) \mathbf{u}_{0, mn}(\tau) - \mathbf{A}_{1, mi}(\tau) \mathbf{u}_{0, mn}(\tau),
\end{equation}
where $\delta_{ij}$ is the Kronecker delta function, 
\begin{align}
\Phi'_{mn}(\tau) &= \sum_{i=1}^{\infty} d_i(\tau) \Phi_{mi}(\tau),\label{Phiseries}\\
\mathbf{A}_1(\tau) \mathbf{u}_{0, mn}(\tau) \Phi_{mn}(\tau) &= \sum_{i=1}^{\infty}\mathbf{A}_{1, mi}(\tau)\mathbf{u}_{0, mn}(\tau) \Phi_{mi}(\tau).\label{u1eq}
\end{align}
With $\Psi'_0(\tau) = \sigma^{\pm}_{mn}(\tau)$, there will be two different cases whether $i = n$ or not. In the latter case $\mathbf{A}_{0, mi}(\tau) - \Psi'_0(\tau) \mathbf{I}$ is non-singular and we may solve (\ref{u1syseq}) for $\mathbf{u}_{1, mi}$. In the case where $i = n$, $\mathbf{A}_0(\tau) - \Psi'_0(\tau) \mathbf{I}$ is singular. We can still solve (\ref{u1syseq}) if the right side of (\ref{u1eq}) satisfies the solvability condition of being orthogonal to any null eigenvector $\mathbf{u}_0^{(*)}(\tau)$ of the transpose of the matrix $\mathbf{A}_{0,mn}(\tau) - \Psi'_0(\tau) \mathbf{I}$. We obtain the following expression 
\begin{equation}\label{dotpsi1}
\Psi'_1(\tau) = \frac{\mathbf{u}_0^{(*)}(\tau)^{\intercal} [ \mathbf{u}'_0(\tau) - d_n(\tau) \mathbf{u}_{0,mn}(\tau) + \mathbf{A}_{1, mn}(\tau) \mathbf{u}_{0,mn}(\tau) ]}{\mathbf{u}_0^{(*)}(\tau)^{\intercal} \mathbf{u}_{0,mn}(\tau)}.
\end{equation}

\subsection{The centre subspace}

We choose parameter values in our model (\ref{GroX}) - (\ref{XYBC}) so that (as in Section \ref{subsec:MargStabCurv}) for the curvature parameter $\gamma = \gamma_0$, the eigenvalues of the constant domain autonomous linearization (\ref{linearization}) are all real and negative, except for one isolated zero eigenvalue $\sigma_{m_0 n_0}^{\pm} = 0$. Furthermore for $\gamma \in [\gamma_1, \gamma_2]$, with $\gamma_1 < \gamma < \gamma_2$, the eigenvalue $\sigma_{m_0 n_0}^{\pm}$ is real and near zero with $\sigma_{m_0 n_0}^{\pm} < 0$ for $\gamma = \gamma_2$ and $\sigma_{m_0 n_0}^{\pm} > 0$ for $\gamma = \gamma_1$, and the remaining eigenvalues are real and negative, uniformly bounded away from zero.

Then we take $\gamma = \gamma(\tau)$ in the slow timescale, decreasing from $\gamma_2$, through $\gamma_0$, to $\gamma_1$. For example, we may take
\begin{equation}\label{gammaepsilon}
\gamma(\tau) = \gamma_0 - \tau,
\end{equation}
for $\tau$ belonging to a suitable interval.

For parameter values corresponding to Figure \ref{fig:gamvsA}, the parameter point $(A, \gamma(\tau))$ moves, as $\tau$ increases, to the left along the dotted arrow that crosses the $(m_0, n_0) = (5, 1)$ marginal stability curve from the regime of stability into the region of instability.

For the slowly changing domain, the leading order $\tau$-dependent linearization (\ref{linearizationtau}) has eigenvalues ${\Psi_0}'(\tau) = \sigma^{\pm}_{mn}(\tau)$ with
\begin{enumerate}
\item $\sigma^{(0)}_0(\tau) = \sigma_{m_0 n_0}^{+}(\tau)$ is the only such value of all the $\sigma^{\pm}_{mn}(\tau)$ that is near $0$, its absolute value remains small, for all $\tau$, is strictly increasing with respect to $\tau$ and goes from negative to positive as $\tau$ is increased.

\item $\sigma_{m_0 n_0}^{-}(\tau)$ and all other $\sigma^{\pm}_{mn}(\tau)$, $(m,n) \neq (m_0, n_0)$ are strictly negative and larger in absolute value than $|\sigma^{(0)}_0(\tau)|$ for all $\tau$.

\end{enumerate}

Then from the WKB approximation (\ref{U0series}), (\ref{U0Phi}), for small $\varepsilon$ the linearization (\ref{WKB}) about the quasi-patternless solution makes a transition, as the slow time $\tau$ increases, from the stable to unstable regimes. We assume for small $\varepsilon$ that the solution space of (\ref{WKB}) splits into a direct sum of two invariant subspaces: a finite dimensional centre subspace associated with $\sigma^{(0)}_0(\tau)$, and a complementary infinite dimensional stable subspace. Solutions of (\ref{WKB}) in the centre subspace decay or grow at a slow rate, while solutions in the stable subspace decay much faster, exponentially, throughout the time interval.

There are two linearly independent critical solutions to the linearized system (\ref{WKBlin}). One of them is
\begin{align}\label{LinPsi}
\mathbf{\hat{U}}^{(0)}(\tau, \varepsilon) = e^{\Psi^{(0)}(\tau, \varepsilon)/\varepsilon} \ \mathbf{U}^{(0)}(\tau, \varepsilon),
\end{align}
where $\frac{\partial}{\partial \tau}\Psi^{(0)}(\tau, \varepsilon) = \sigma_0^{(0)}(\tau) + \mathcal{O}(\varepsilon)$, and the other critical solution is the complex conjugate of (\ref{LinPsi}).

Let $\mathcal{X}^c(\tau, \varepsilon)$ be the span of $\hat{\mathbf{U}}^{(0)}(\tau, \varepsilon)$ and its complex conjugate, then this space is two-dimensional and may be expressed as
\begin{align}\label{centrespace}
\mathcal{X}^c(\tau, \varepsilon) &= \left\{ z \mathbf{U}^{(0)}(\tau, \varepsilon) + \bar{z} \bar{\mathbf{U}}^{(0)}(\tau, \varepsilon) | z \in \mathbb{C} \right\}.
\end{align}
The action of the symmetries (\ref{Symm}) on the space $\mathcal{X}^c(\tau, \varepsilon)$ has the following effects on $z$:
\begin{align}
z &\to e^{i \varphi_0} z, & z &\to \bar{z}.\label{zsyms}
\end{align}



Let $\mathbf{P}^c(\tau, \varepsilon)$ be a projection onto $\mathcal{X}^c(\tau, \varepsilon)$, which we assume has the form
\begin{align}
\mathbf{P}^c(\tau, \varepsilon) = \mathbf{P}_0^c(\tau) + \varepsilon \mathbf{P}_1^c(\tau) + \mathcal{O}(\varepsilon^2).
\end{align}
We will only need the leading order term $\mathbf{P}_0^c(\tau)$, given by
\begin{align}\label{Projectioneq}
\mathbf{P}_0^c(\tau) = \langle \mathbf{U}^{(*)}_0(\tau), \mathbf{U} \rangle \mathbf{U}^{(0)}_0(\tau) + \langle \bar{\mathbf{U}}^{(*)}_0(\tau), \mathbf{U} \rangle \bar{\mathbf{U}}^{(0)}_0(\tau),
\end{align}
where the angled brackets denote the inner product
\begin{align}
\begin{split}
\langle \mathbf{U}_1, \mathbf{U}_2 \rangle &= \int_0^{2\pi} \int_0^{\arcsin \gamma (\tau)} (\bar{U}_1 U_2 + \bar{V}_1 V_2) \sin \theta \, d \theta \, d \varphi\\
&= \int_0^{2\pi} \int_0^{\infty} (\bar{U}_1 U_2 + \bar{V}_1 V_2) \frac{\sinh \eta \sin^2 \xi(\tau)}{[\cosh \eta - \cos \xi(\tau)]^2} \, d \eta \, d \varphi,
\end{split}
\end{align}
for $\mathbf{U}_j = (U_j, V_j)^{\intercal}$, $j = 1, 2$, and
\begin{align}
\mathbf{U}^{(*)}_0(\tau) = N^{(*)}(\tau) \begin{pmatrix} k_3 \\ D_X \mu_{m_0 n_0}(\tau) - k_1 + \sigma^{(0)}_0(\tau) \end{pmatrix} \Phi_{m_0 n_0}(\tau)
\end{align}
is the solution of the adjoint problem to (\ref{linearizationtau}) with $m=m_0$, $n=n_0$, $\Psi'_0(\tau) = \sigma^{(0)}_0(\tau)$ and the normalization coefficient $N^{(*)}(\tau)$ selected so that
\begin{align}
\langle \mathbf{U}_0^{(*)}(\tau), \mathbf{U}_0^{(0)}(\tau) \rangle = 1
\end{align}
for all $\tau$. Let $\mathcal{X}^s(\tau, \varepsilon)$ be the kernel of $P^c(\varepsilon, \tau)$, then the projection onto $\mathcal{X}^s(\tau, \varepsilon)$ is
\begin{align}
\mathbf{P}^s(\tau, \varepsilon) = \mathbf{I} - \mathbf{P}^c(\tau, \varepsilon).
\end{align}


\subsection{Centre manifold equations}

Following the analysis in the constant domain case, for the slowly changing domain we reduce the nonlinear partial differential equation system (\ref{ABCeq}) into an ordinary differential equation on a centre manifold. The solutions of the reduced system give a good approximation of the transition to the patterned state. In this case, however, the equations are nonautonomous, so the normal form coefficients we obtain after reduction will depend on time. 

We split the function space $\mathcal{X}$ of all sufficiently regular functions on the spherical cap satisfying homogeneous Dirichlet boundary conditions into the two subspaces
\begin{align}\label{spacesplit}
\mathcal{X} &= \mathcal{X}^c(\tau, \varepsilon) \oplus \mathcal{X}^s(\tau, \varepsilon).
\end{align}
We assume the splitting (\ref{spacesplit}) is invariant under the linearization (\ref{ABCeq}) for all $\tau$, $\varepsilon$. Each function $\mathbf{U} \in \mathcal{X}$ may then be written uniquely as a sum of two components belonging to each subspace and each component may be retrieved using projections,
\begin{align}
\mathbf{U} = \mathbf{U}^c(\tau, \varepsilon) + \mathbf{U}^s(\tau, \varepsilon),
\end{align}
where
\begin{align}
\mathbf{U}^c(\tau, \varepsilon) &= \mathbf{P}^c(\tau, \varepsilon)\mathbf{U} \in \mathcal{X}^c(\tau, \varepsilon), & \mathbf{U}^s(\tau, \varepsilon) &= \mathbf{P}^s(\tau, \varepsilon)\mathbf{U} \in \mathcal{X}^s(\tau, \varepsilon).
\end{align}
Thus, equation (\ref{ABCeq}) may be rewritten as the system
\begin{align}
\frac{\partial\mathbf{U}^c}{\partial t} &= \mathbf{A}(\varepsilon t, \varepsilon) \mathbf{U}^c + \mathbf{P}^c(\varepsilon t, \varepsilon) \mathbf{F}(\mathbf{U}^c(\varepsilon t, \varepsilon) + \mathbf{U}^s(\varepsilon t, \varepsilon), \varepsilon t, \varepsilon)\label{UcABCeq}\\
\frac{\partial\mathbf{U}^s}{\partial t} &= \mathbf{A}(\varepsilon t, \varepsilon) \mathbf{U}^s + \mathbf{P}^s(\varepsilon t, \varepsilon) \mathbf{F}(\mathbf{U}^c(\varepsilon t, \varepsilon) + \mathbf{U}^s(\varepsilon t, \varepsilon), \varepsilon t, \varepsilon),\label{UsABCeq}
\end{align}
with $\mathbf{F}(\mathbf{U}, \varepsilon t, \varepsilon) = \mathbf{B}(\varepsilon t, \varepsilon)(\mathbf{U}, \mathbf{U}) + \mathbf{C}(\mathbf{U}, \mathbf{U}, \mathbf{U})$. We assume that there is an exponentially attracting centre manifold tangent to $\mathcal{X}^c(\varepsilon t, \varepsilon)$
\begin{align} \label{EACM}
\mathbf{U}^s &= \mathbf{W}(\mathbf{U}^c, \varepsilon t, \varepsilon)
\end{align}
for all $\mathbf{U}^c$ in the subspace $\mathcal{X}^c(\varepsilon t, \varepsilon)$ with sufficiently small magnitude. Equation (\ref{UsABCeq}) is then written as (cf. \cite{Potzsche2006})
\begin{align}
\begin{split}
\frac{\partial}{\partial t} \mathbf{W}(\mathbf{U}^c, \varepsilon t, \varepsilon) &+ \frac{\partial}{\partial \mathbf{U}^c} \mathbf{W}(\mathbf{U}^c, \varepsilon t, \varepsilon) \left[ \mathbf{A}(\varepsilon t, \varepsilon) \mathbf{U}^c + \mathbf{P}^c(\varepsilon t, \varepsilon) \mathbf{F}(\mathbf{U}^c + \mathbf{W}(\mathbf{U}^c, \varepsilon t, \varepsilon), \varepsilon t, \varepsilon) \right]\\
&= \mathbf{A}(\varepsilon t, \varepsilon) \mathbf{W}(\mathbf{U}^c, \varepsilon t, \varepsilon)
+ \mathbf{P}^s(\varepsilon t, \varepsilon) \mathbf{F}(\mathbf{U}^c + \mathbf{W}(\mathbf{U}^c, \varepsilon t, \varepsilon), \varepsilon t, \varepsilon),
\end{split}\label{dWdt1}
\end{align}
after substituting for $\frac{\partial \mathbf{U}^c}{\partial t}$. The Taylor expansion of $\mathbf{W}$ around $\mathbf{U}^c = 0$ starts at the second degree due to the tangency condition
\begin{align}
\mathbf{W}(\mathbf{U}^c, \varepsilon t, \varepsilon) = \mathbf{W}^{(1)}(\varepsilon t, \varepsilon) (\mathbf{U}^c, \mathbf{U}^c) + \mathcal{O}(\| \mathbf{U}^c \|^3),
\end{align}
where $\mathbf{W}^{(1)}(\varepsilon t, \varepsilon)$ is a symmetric bilinear form. Applying the product rule to the $\mathbf{U}^c$ derivative yields
\begin{align}
\frac{\partial}{\partial \mathbf{U}^c} \mathbf{W}(\mathbf{U}^c, \varepsilon t, \varepsilon) \mathbf{V} = 2 \mathbf{W}^{(1)}(\varepsilon t, \varepsilon)(\mathbf{U}^c, \mathbf{U}^c) + O(\|\mathbf{U}^c\|\|\mathbf{V}\|)
\end{align}
for any vector $\mathbf{V}$. 

As a result of equivariance of the system of differential equations (\ref{ConUvec}) under rotations and reflections (\ref{zsyms}) we may only study the real part of $z$, labeled $x$ and then apply a phase $\phi$ to get the other solutions $z = xe^{i\phi}$. Due to this we will continue with using $x$ in our derivations as a solution representative of the family of solutions.

If we choose $\mathbf{U}^c = x \tilde{{\mathbf{U}}}^{(0)}$, where $\tilde{\mathbf{U}}^{(0)} = \text{Re } \mathbf{U}^{(0)}$ for $x$ real then equation (\ref{dWdt1}) can be used to find the order $x^2$ terms with
\begin{align}
\begin{split}
\frac{\partial}{\partial t} \mathbf{W}^{(1)}&(\varepsilon t, \varepsilon) (\mathbf{U}^{(0)}(\varepsilon t, \varepsilon), \mathbf{U}^{(0)}(\varepsilon t, \varepsilon)) + 2\mathbf{W}^{(1)} (\mathbf{U}^{(0)}(\varepsilon t, \varepsilon), \mathbf{A}(\varepsilon t, \varepsilon) \mathbf{U}^{(0)}(\varepsilon t, \varepsilon)) \\
&= \mathbf{A}(\varepsilon t, \varepsilon) \mathbf{W}^{(1)}(\varepsilon t, \varepsilon) (\mathbf{U}^{(0)}(\varepsilon t, \varepsilon), \mathbf{U}^{(0)}(\varepsilon t, \varepsilon)) + \mathbf{P}^s(\varepsilon t, \varepsilon) \mathbf{B}(\mathbf{U}^0(\varepsilon t, \varepsilon), \mathbf{U}^{(0)}(\varepsilon t, \varepsilon)),
\end{split}\label{dWdt2}
\end{align}
dropping the tilde for the real part of $\mathbf{U}^{(0)}$. By defining
\begin{align}
\mathbf{U}^{(1)}(\varepsilon t, \varepsilon) = \mathbf{W}^{(1)}(\varepsilon t, \varepsilon)(\mathbf{U}^{(0)}(\varepsilon t, \varepsilon), \mathbf{U}^{(0)}(\varepsilon t, \varepsilon)),
\end{align}
we can write the solution $\mathbf{U}$ as a series in $x$
\begin{align}\label{centremanifold}
\mathbf{U} = x \mathbf{U}^{(0)}(\varepsilon t, \varepsilon) + x^2 \mathbf{U}^{(1)}(\varepsilon t, \varepsilon) + O(x^3).
\end{align}
We use $\varepsilon t = \tau$ as well as
\begin{align}
\mathbf{A}(\tau, \varepsilon) \mathbf{U}^{(0)}(\tau, \varepsilon) = \frac{\partial \Psi^{(0)}}{\partial \tau}(\tau, \varepsilon)\mathbf{U}^{(0)}(\tau, \varepsilon) + \varepsilon\frac{\partial}{\partial \tau}\mathbf{U}^{(0)}(\tau, \varepsilon) \label{AU0again}
\end{align}
and
\begin{align}
\frac{\partial}{\partial \tau}\mathbf{U}^{(1)}(\tau, \varepsilon) = \frac{\partial}{\partial \tau}\mathbf{W}^{(1)}(\tau, \varepsilon) (\mathbf{U}^{(0)}(\tau, \varepsilon), \mathbf{U}^{(0)}(\tau, \varepsilon)) + 2\mathbf{W}^{(1)}(\tau, \varepsilon)  \left( \mathbf{U}^{(0)}(\tau, \varepsilon), \frac{\partial}{\partial \tau} \mathbf{U}^{(0)}(\tau, \varepsilon) \right),
\end{align}
to get the expression
\begin{align}
\left[ \varepsilon \frac{\partial}{\partial \tau} - \mathbf{A}(\tau, \varepsilon) + 2\frac{\partial \Psi^{(0)}}{\partial \tau}(\tau, \varepsilon) \right]\mathbf{U}^{(1)}(\tau, \varepsilon) = \mathbf{P}^s(\tau, \varepsilon) \mathbf{B}(\tau, \varepsilon)(\mathbf{U}^{(0)}(\tau, \varepsilon), \mathbf{U}^{(0)}(\tau, \varepsilon)),
\end{align}
which we can use to solve for $\mathbf{U}^{(1)}(\varepsilon t, \varepsilon) = \mathbf{U}^{(1)}_0(\varepsilon t) + \mathcal{O}(\varepsilon)$, assuming that $\mathbf{U}^{(1)}_0(\varepsilon t)$ has a series expression
\begin{align}
\mathbf{U}^{(1)}_0(\varepsilon t) = \sum_{j=1}^{\infty} \left[\mathbf{u}^{(1)}_{0j}(\varepsilon t) \mathbf{\Phi}_{0j}(\varepsilon t) + \mathbf{u}^{(1)}_{2mj}(\varepsilon t) \mathbf{\Phi}_{2mj}(\varepsilon t) \right].\label{U1series}
\end{align}

\subsection{Normal form}\label{subsec:NF}

Having found the leading terms in the centre manifold, we can use a projection on the centre subspace $\mathcal{X}^c(\varepsilon t, \varepsilon)$. This is given by equation (\ref{UcABCeq})
\begin{align}
\frac{\partial}{\partial t} \mathbf{U}^{c} + \mathbf{P}^c(\varepsilon t, \varepsilon) \mathbf{F}(\mathbf{U}^c + \mathbf{U}^s, \varepsilon t, \varepsilon).\label{UsABCeq2}
\end{align}
We proceed here as if $\mathbf{W}^c(\varepsilon t, \varepsilon)$ is such that
\begin{align}
\mathbf{U}^c &= x \mathbf{U}^{(0)}(\varepsilon t, \varepsilon), & \mathbf{U}^s &= x^2 \mathbf{U}^{(1)}(\varepsilon t, \varepsilon) + \mathcal{O}(x^3).
\end{align}
Using this in (\ref{UsABCeq2}) yields
\begin{align}
\begin{split}
\frac{dx}{dt} \mathbf{U}^{(0)}(\varepsilon t, \varepsilon) + x \varepsilon \frac{\partial}{\partial \tau}&\mathbf{U}^{(0)}(\varepsilon t, \varepsilon) = x \mathbf{A}(\varepsilon t, \varepsilon) \mathbf{U}^{(0)}(\varepsilon t, \varepsilon)\\
 &+ \mathbf{P}^c(\varepsilon t, \varepsilon) \mathbf{F}(x \mathbf{U}^{(0)}(\varepsilon t, \varepsilon) + x^2 \mathbf{U}^{(1)}(\varepsilon t, \varepsilon) + \mathcal{O}(x^3), \varepsilon t, \varepsilon).
\end{split}
\end{align}
We then use (\ref{AU0again}) and carry out a projection in the direction of $\mathbf{U}^{(0)}(\varepsilon t, \varepsilon)$ to obtain
\begin{align}\label{1DNF} 
\dot{x} &= \sigma^{(0)}(\varepsilon t, \varepsilon) x + C(\varepsilon t, \varepsilon) x^3 + \mathcal{O}(x^3),
\end{align}
with
\begin{align}
\sigma^{(0)}(\varepsilon t, \varepsilon) &= \frac{\partial}{\partial \tau} \Psi^{(0)} (\varepsilon t, \varepsilon) = \sigma^{(0)}_0(\varepsilon t) + \varepsilon \sigma^{(0)}_1(\varepsilon t) + \mathcal{O}(\varepsilon^3),\\
C(\varepsilon t, \varepsilon) &= C_0(\varepsilon t) + \mathcal{O}(\varepsilon)\\
C_0(\varepsilon t) &= \left\langle \mathbf{U}^{(*)}_0(\varepsilon t), 2\mathbf{B}_0 \left( \mathbf{U}^{(0)}_0(\varepsilon t), \mathbf{U}^{(1)}_0(\varepsilon t) \right) + \mathbf{C}_0 \left( \mathbf{U}^{(0)}_0(\varepsilon t), \mathbf{U}^{(0)}_0(\varepsilon t), \mathbf{U}^{(0)}_0(\varepsilon t) \right) \right\rangle.\label{C0eqfin}
\end{align}

Equations (\ref{linearizationtau}), (\ref{dotpsi1}) and (\ref{C0eqfin}) allow us to find approximate normal form coefficients for any time $t$. In order to use these equations we need a good approximation for the infinite sums of eigenfunctions of modes $0$, $m_0$ or $2m_0$ present, for example, in the quadratic component of the centre manifold $\mathbf{U}^{(1)}$ or the derivative of the eigenfunction $\mathbf{\Phi}'_{m_0n}(\tau)$, both of which are found through equations (\ref{U1series}) and (\ref{Phiseries}). For $m_0 = 5$ we found that truncating the series after the first five terms for each order yield similar results to computing the first eight or twelve terms, which suggests that truncated series approximations are sufficiently accurate. 

\subsubsection{Solving the normal form equations}
To solve the normal form ordinary differential equation, we need the coefficient values throughout the time interval. As the coefficients do not have an algebraic form that can be easily solved using symbolic computing tools, we cannot solve at all time values. To get a good approximation, we collect the coefficient values at a sample of equally spaced times and interpolate between them. We used a sample of 36 time points and spline interpolation using the \emph{CurveFitting} tool in the Maple computational software package, which uses piecewise cubic polynomials to create an interpolation that matches all sample point values and is continuous up its second derivative. 

\begin{figure}
	\centering
		\includegraphics[width=1.00\textwidth]{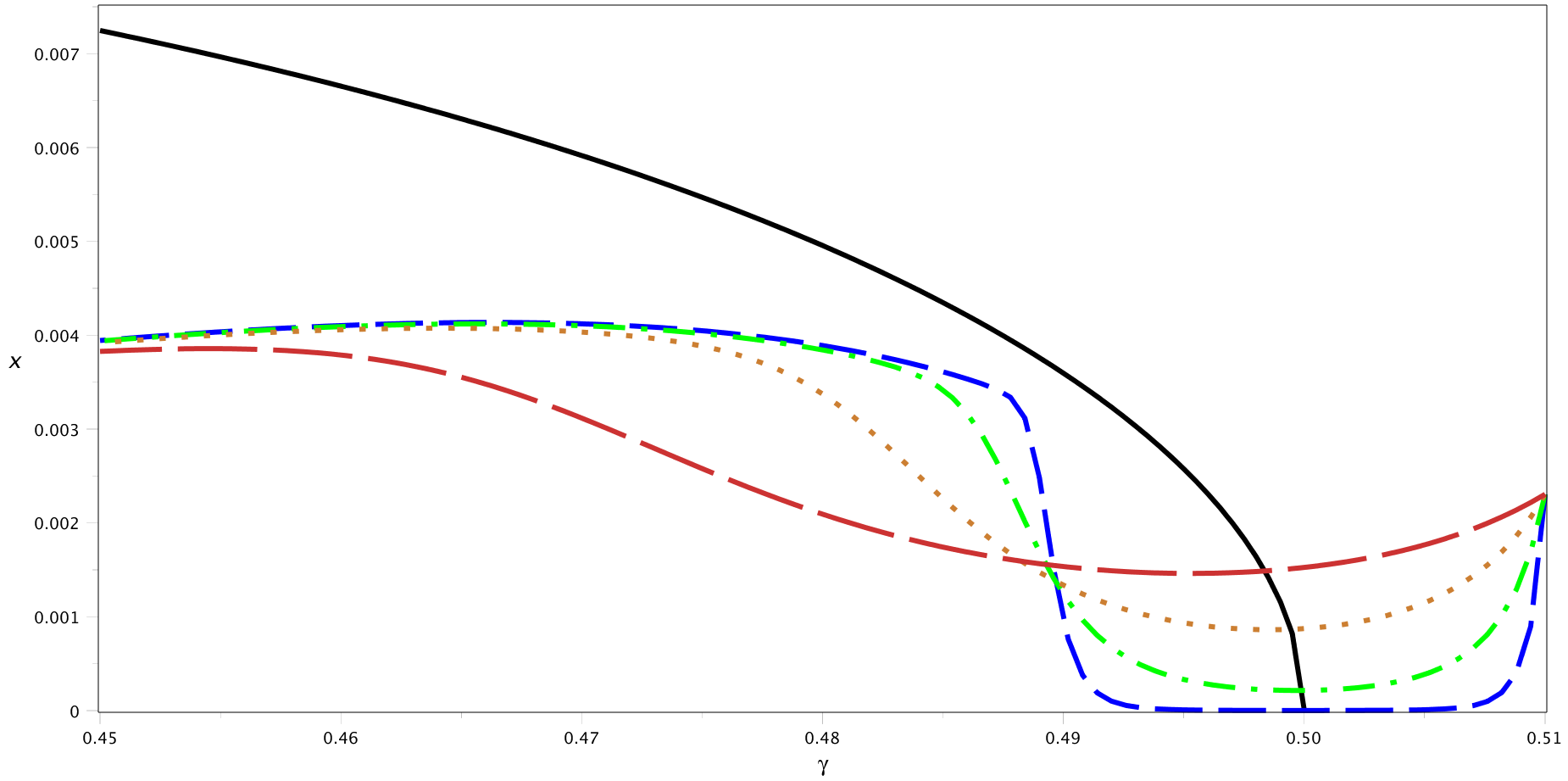}
	\caption{Normal form solutions and the constant domain bifurcation diagram for parameter values (\ref{param1}), (\ref{param2}), $A = 76.5198$ and different $\varepsilon$ values. The values of $\varepsilon$ are $1 \times 10^{-6}$ for the long-dashed orange curve, $3 \times 10^{-7}$ for the dotted gold curve, $1 \times 10^{-7}$ for the dot-dashed green curve and $3 \times 10^{-8}$ for the dashed blue curve. Each solution has the same initial conditions $\gamma(0) = 0.51$, $x(0) = 0.002305$. The solid black curve is the constant coefficient pitchfork bifurcation diagram.}
	\label{fig:NFepsplotssame}
\end{figure}

Figure \ref{fig:NFepsplotssame} show solutions to the normal form equations for different $\varepsilon$ values superimposed with the constant domain pitchfork bifurcation nontrivial stable equilibrium branch. We clearly see that a time dependent cubic coefficient $C(\tau, \varepsilon)$ is important to get a good prediction out of our normal form reduction.

The normal form solutions were computed using the \emph{dsolve} function in the Maple software package using a fourth-and-fifth order Runge-Kutta method with adaptive time stepping. We chose the initial condition $x(0) = 2.3057 \cdot 10^{-3}$ to match the initial condition used in the numerical simulations.

\subsection{Numerical simulations}\label{sec:Numerical}

To validate the normal form predictions we computed numerical simulations of solutions to the nonautonomous Brusselator system (\ref{RDEvec}) -- (\ref{RDEBC}) using the closest point method \cite{RuuthMerriman2008} with implicit-explicit time stepping \cite{Ascher1995}. We used parameter values (\ref{param1}), (\ref{param2}) and
\begin{align}
A &= 76.51981.
\end{align} 

\begin{figure}[ht]
	\centering
	\subfloat{\label{fig:QP_dx_02_top}\includegraphics[width=0.5\linewidth, valign=c, trim={0cm 0 0cm 0},clip]{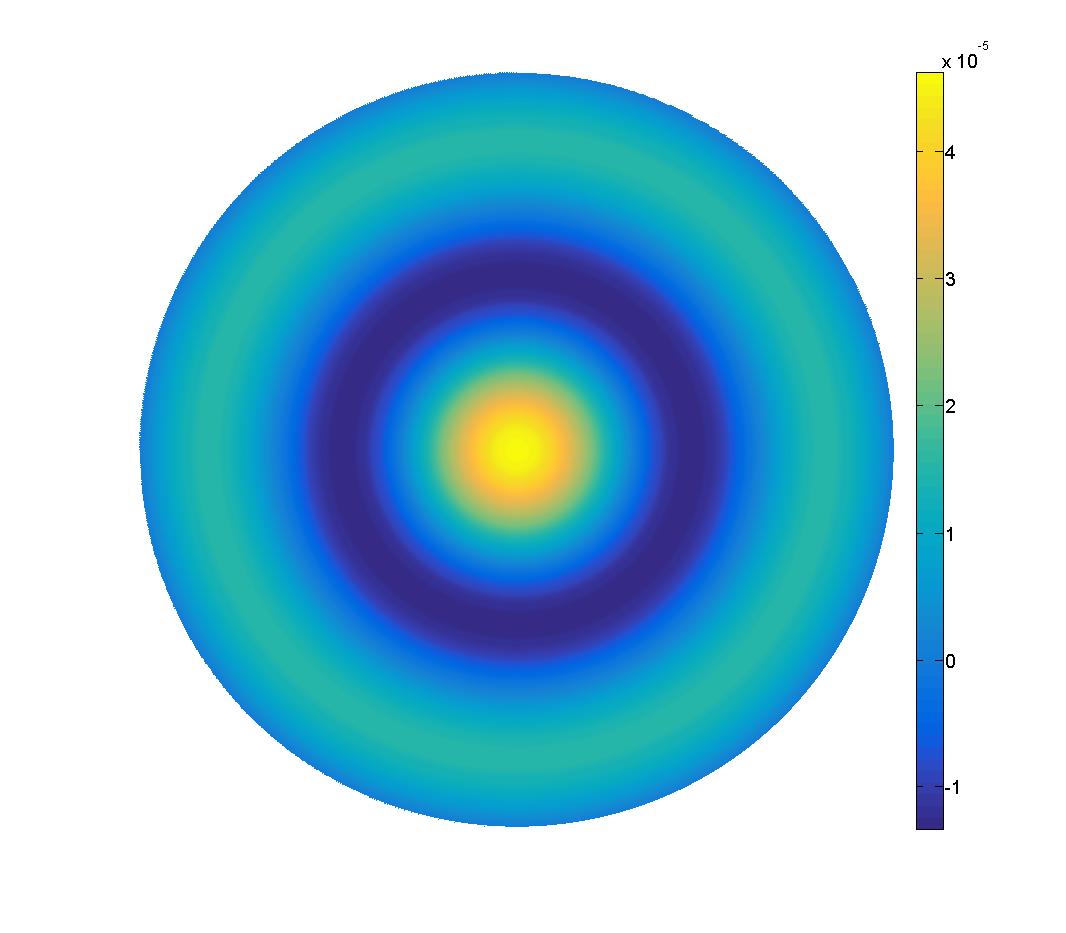}}
	\subfloat{\label{fig:QP_dx_02_3D}\includegraphics[width=0.5\linewidth, valign=c, trim={0cm 0 0cm 0},clip]{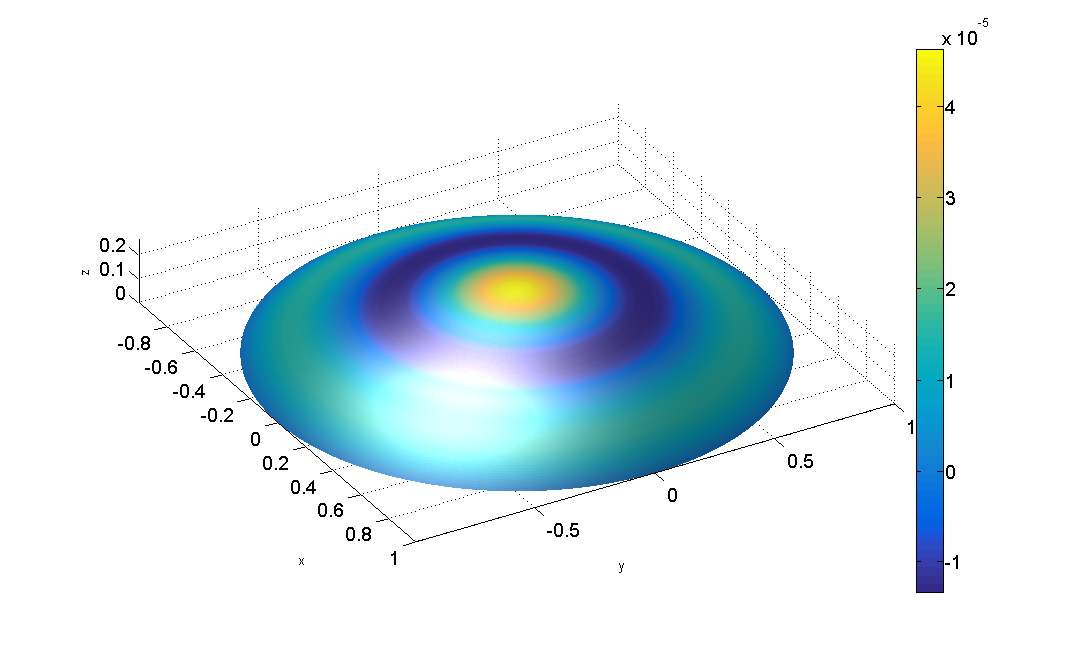}}
	\caption{Quasi-patternless solution of (\ref{ConUvec}) with modified affine operator (\ref{QPAmod}), simulated using the closest point method, after $10000$ time units from a homogeneous initial condition with $h = 0.02$, $\Delta t = 0.1$. Parameter values are $\varepsilon = 10^{-6}$, $\gamma = 0.4915$, $A = 76.51981$, and the same as (\ref{param1}), and (\ref{param2}).}
	\label{fig:QP_dx_02}
\end{figure}

\begin{figure}[ht]
	\centering
		\includegraphics[width=0.80\textwidth]{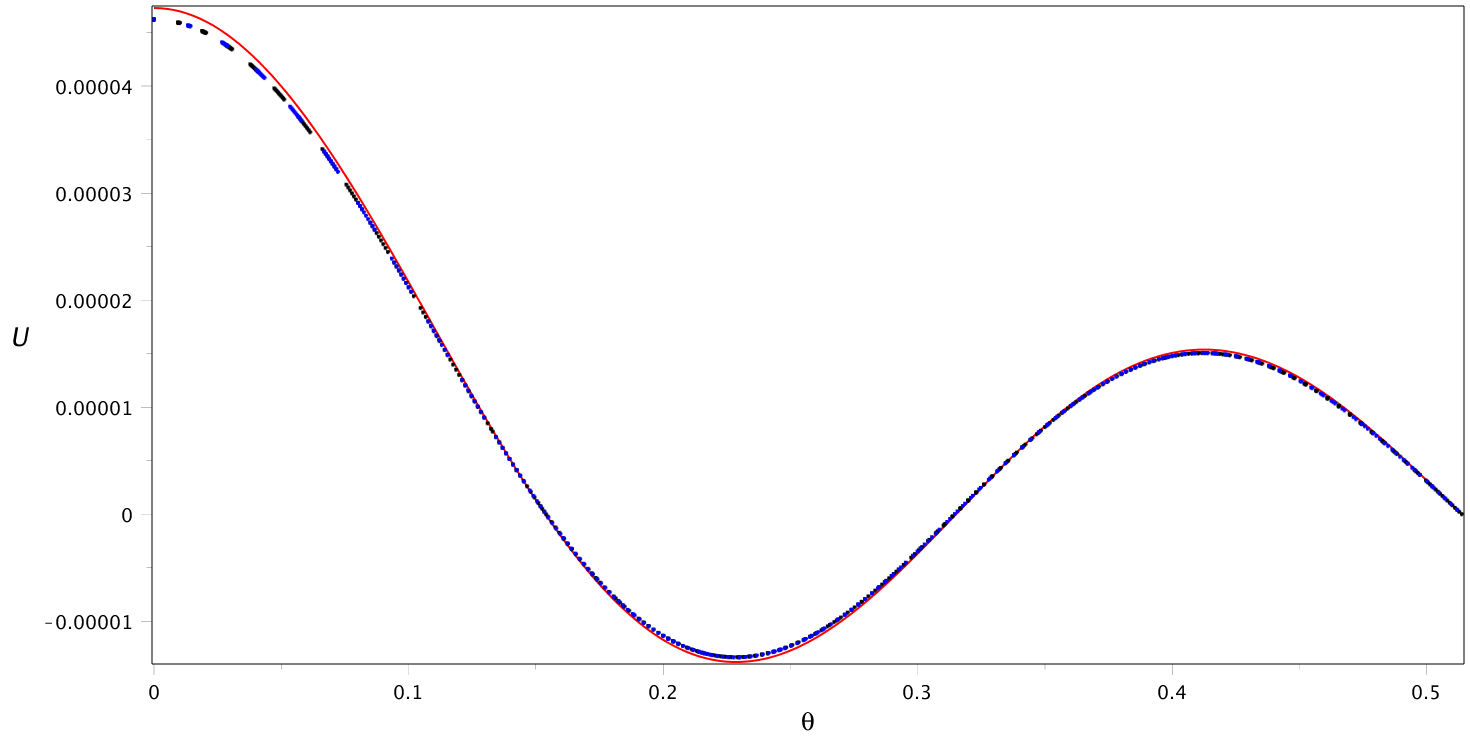}
	\caption{Quasi-patternless solution simulation computed using the closest point method (blue dots) for $h = 0.02$ with the same parameter values as in Figure \ref{fig:QP_dx_02}, compared to a five-term truncated series approximation (red curve) outlined in (\ref{QPseries}) with coefficients as computed in (\ref{alphalinsol}).}
	\label{fig:QP_02_diag}
\end{figure}

We first checked the accuracy of our truncated series approximations of the quasi-patternless solution, by computing simulations of a modified version of (\ref{ABCeq}), where the linear operator $\mathbf{A}$ is replaced with the affine linear operator defined by
\begin{align}\label{QPAmod}
\tilde{\mathbf{A}}(\tau, \varepsilon)\mathbf{U} = \mathbf{D} \Delta_{\Omega(\tau)} \mathbf{U} + \mathbf{K}(\tau, \varepsilon)\mathbf{U} - \varepsilon \gamma'(\tau) Q(\tau)\left[\mathbf{U} + \mathbf{X}_{00}\right], 
\end{align}
where $\mathbf{U}$ is the deviation from the patternless solution of the constant domain system. We used the same parameter values except that we had our curvature index $\gamma(\tau)$ go from $0.5015$ to $0.4915$, where the quasi-patternless solution remains stable. The simulations ran for $10000$ time units, with $\varepsilon = 10^{-6}$, and with a spatial grid size $h = 0.02$ and time step size of $0.1$. We extracted the numerical solution of $\mathbf{U}$ over the polar angle $\theta$ for two longitudinal angles $\varphi$ separated by $\pi/4$ and compared it with the five-term truncated series from our normal form calculations. The truncated series agrees well with the simulations and the simulations converge with decreasing spatial step $h$. Truncated series with twelve terms were also compared and yielded very similar results. Figures \ref{fig:QP_dx_02} and \ref{fig:QP_02_diag} show a numerical simulation that show the quasi-patternless solution on a spherical cap and a cross-section along one longitudinal angle respectively. In Figure \ref{fig:QP_02_diag} the simulation is compared to a truncated series approximation.

We computed other numerical simulations of the full system (\ref{ABCeq}) that show the delayed bifurcation. We ran a typical simulation for $40000$ time units, with $\varepsilon = 10^{-6}$, and with a time step size of $0.1$ and different grid sizes between $0.04$ and $0.015$. Due to the slow domain evolution, we updated the domain geometry only every 50 time steps instead of every time step to save computation time. We used an initial condition 
\begin{align}\label{PDEIC}
\mathbf{U}(0) = \frac{0.021}{u_{0, 51} \cdot \tilde{M}_{51}} \cdot \mathbf{u}_{0, 51} \cdot \tilde{\Phi}_{51},
\end{align}
where $\tilde{\Phi}_{51}$ is the approximate $(5, 1)$ eigenfunction of the Laplace-Beltrami operator computed using the closest point method \cite{Macdonald2011} and the \textsc{MATLAB} \emph{eigs} function, $\tilde{M}_{51}$ is its approximate maximal value, $\mathbf{u}_{0, mn}$ is computed by (\ref{u0eigv}) and $u_{0, 51} = k_2$ is its first coordinate. This makes $0.021$ the maximal value of the initial condition. We use $\varepsilon  = 10^{-6}$ and the curvature starts at $\gamma = 0.4915$ and ends at $\gamma = 0.4515$. We start below the constant domain critical curvature $\gamma_0 = 0.5$ because, as seen in section \ref{subsec:NF} the order $\varepsilon$ terms shift the transition value for $\gamma$, where solutions go from decreasing to increasing in time, to a lower value. We want to start as close as possible to this transition value in order to see solutions ``level off'' to the patterned relative equilibrium at later times, while also keeping computation times reasonable.

To make the comparison between the numerical simulations and normal form solutions, we extract the $(5, 1)$ part of the numerical simulation using the projection (\ref{Projectioneq}) via numerical integration. We then solve for the $x$ value by dividing the maximum of the residue by $u_{0, 51}(\tau)$ times the maximum of the Legendre function at each time. For example consider the initial condition stated in (\ref{PDEIC}) with an exact eigenfunction $\Phi_{51}$. It lies in the centre subspace, therefore the projection will yield the same function
\begin{align}
\mathbf{P}^c \mathbf{U}(0) = \mathbf{U}(0).
\end{align}
Assuming the centre subspace and centre manifold to be close at this value we use the real version of the centre subspace definition (\ref{centrespace}) and comparing
\begin{align}
\mathbf{U}(0) = x(0) \mathbf{U}^{(0)}(0, \varepsilon) = \mathbf{U}_0 = \frac{0.021}{u_{0,51} \cdot M_{51}(0)} \cdot \mathbf{U}^{(0)}(0, \varepsilon),
\end{align}
we obtain
\begin{align}
x(0) = \frac{0.021}{k_2 \cdot M_{51}(0)} \approx 2.3057 \cdot 10^{-3}.
\end{align}
For a later time $t$ a similar computation is used to find $x(t)$,
\begin{align}
x(t) = \| \mathbf{U}(t) \|_{51}{k_2 \cdot M_{51}(\varepsilon t)},
\end{align}
where $\| \cdot \|_{51}$ denotes the $L^{\infty}$ norm of the $(5,1)$-mode projection and $M_{51}(\varepsilon t)$ is the absolute maximum of the $(5,1)$ Legendre function for curvature $\gamma(\varepsilon t)$.

\begin{figure}[ht]
	\centering
		\includegraphics[width=1.00\textwidth]{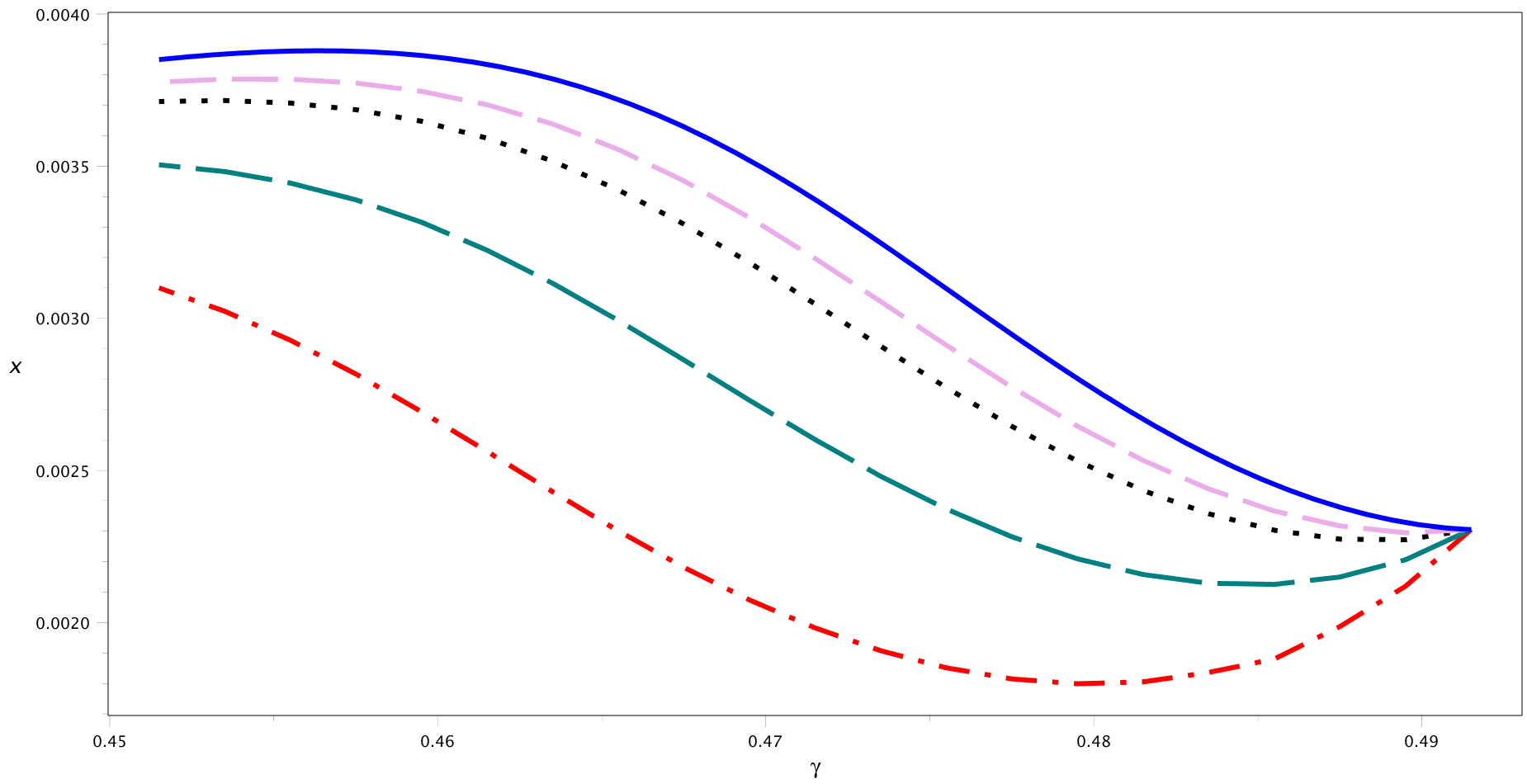}
	\caption{Numerical simulations compared to the normal form solution as a function of curvature index $\gamma$ by extracting $(5,1)$ coefficients using a projection. The spatial grid spacings are $h = 0.04$ for the red dash-dotted curve, $0.03$ for the turquoise long dashed curve, $0.02$ for the black dotted curve and $0.015$ for the purple dashed line and $\varepsilon = 10^{-6}$ for each curve. The normal form solution is the blue solid curve. The horizontal axis represents the curvature index $\gamma$ so solutions go from right to left as time increases. The vertical axis is the variable $x$ for the normal form or the appropriately scaled $L^{\infty}$ norm of the $(5,1)$-mode projection of the numerical solutions.}
	\label{fig:CPM_slns}
\end{figure}

Figure \ref{fig:CPM_slns} shows the resulting equivalent $x(t)$ values for the closest point method simulations for spatial grid sizes between $h = 0.04$ and $h = 0.015$ (dashed or dotted lines) compared with the equivalent normal form computation (solid line). The two initial conditions are chosen so that each numerical simulation match the normal form solution to leading order in $x$ and $\varepsilon$ at the initial time. We see that as the spatial grid is refined the simulations approach the results of the normal form computation. 

\begin{figure}[tbhp]
	\centering
	\subfloat[$\gamma = 0.4915$]{\label{fig:dx_02_rand3_4915}\includegraphics[width=0.5\linewidth]{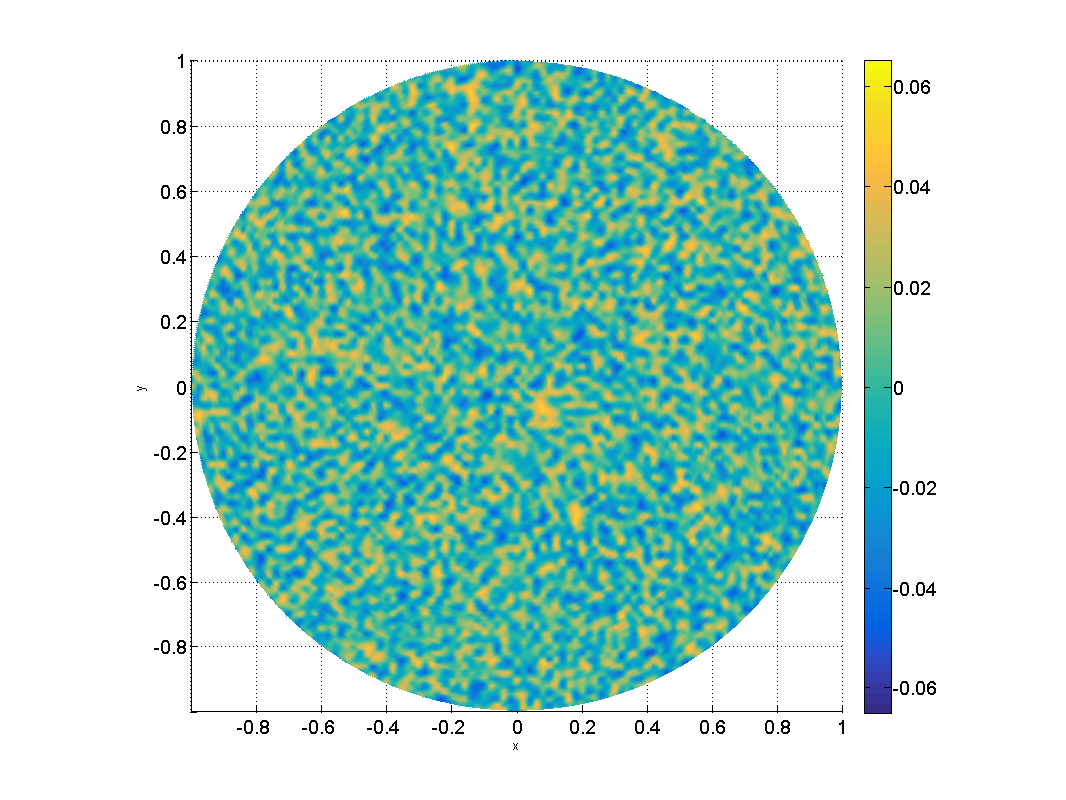}}
	\subfloat[$\gamma = 0.4885$]{\label{fig:dx_02_rand3_4885}\includegraphics[width=0.5\linewidth]{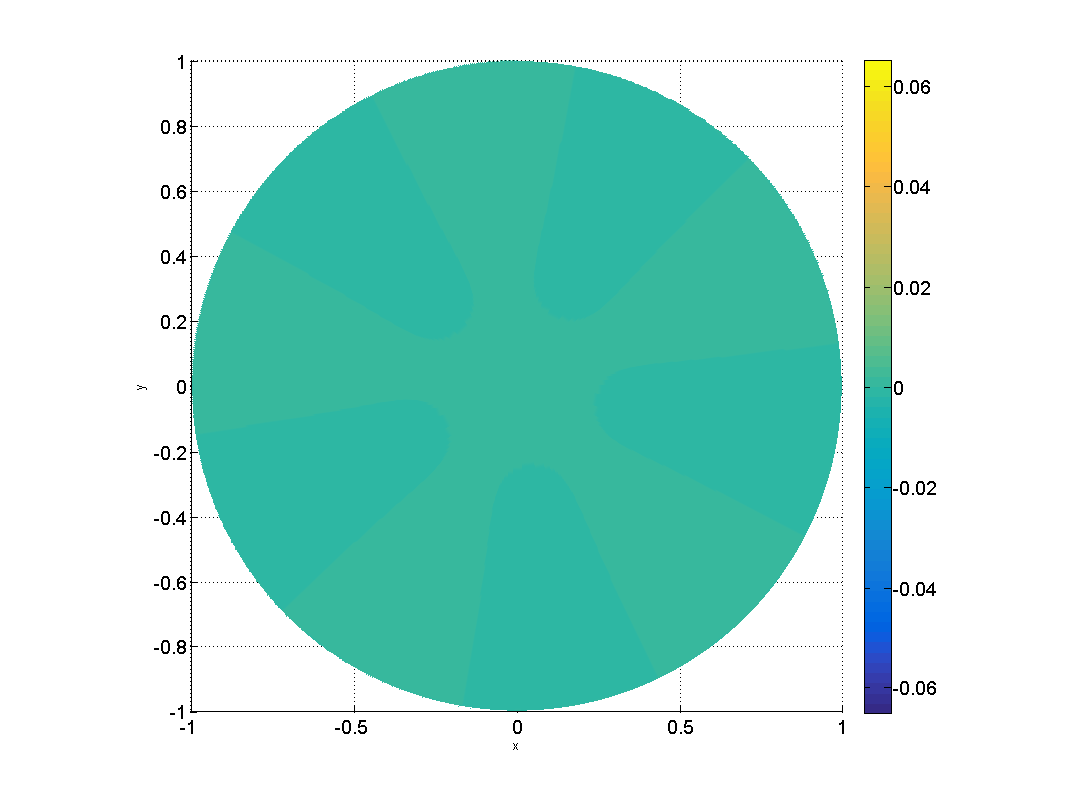}}\\
	\subfloat[$\gamma = 0.4465$]{\label{fig:dx_02_rand3_4465}\includegraphics[width=0.5\linewidth]{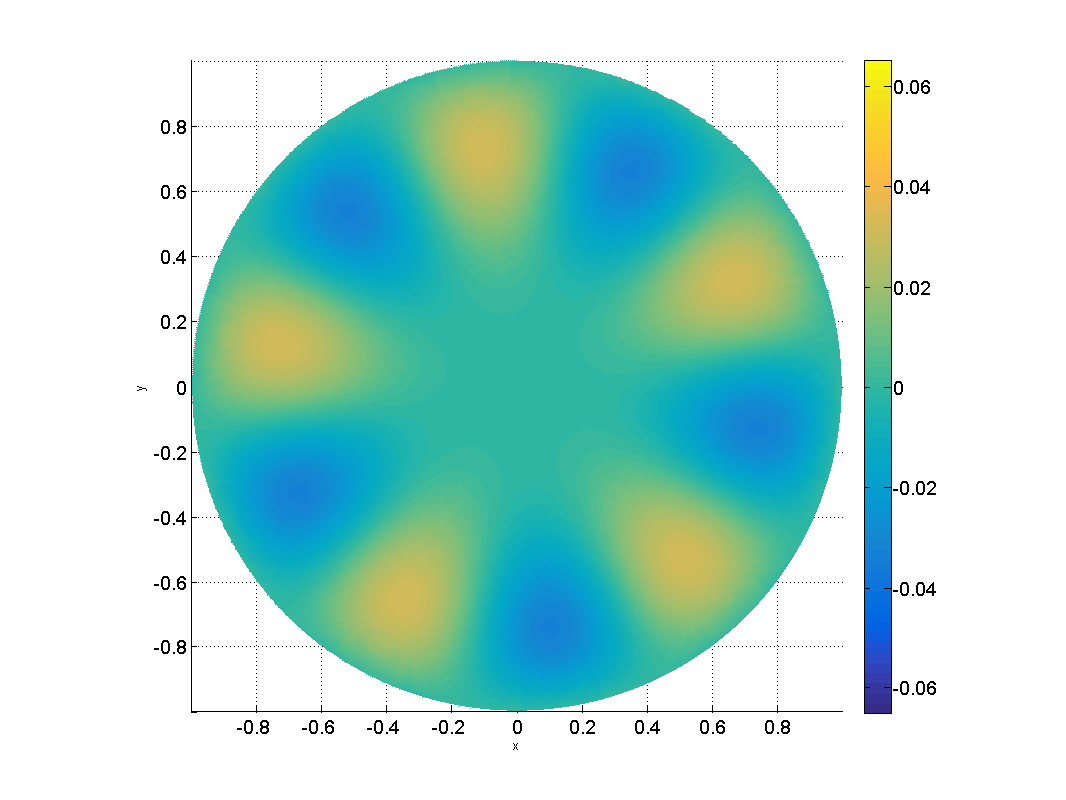}}
	\subfloat[$\gamma = 0.4315$]{\label{fig:dx_02_rand3_4315}\includegraphics[width=0.5\linewidth]{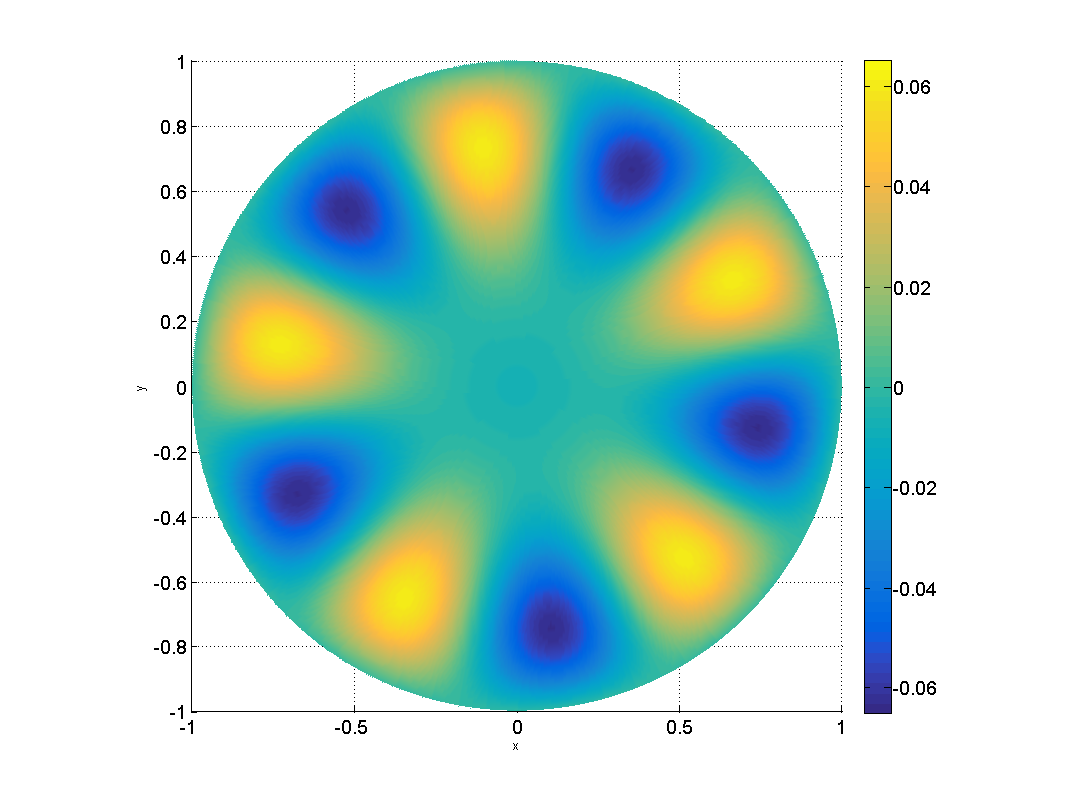}}
	\caption{Snapshots of the $U$-component of numerical simulations using grid size $h = 0.02$ and time step $\Delta t = 0.1$ running for $60000$ time units at four different times. Parameter values used are listed in (\ref{param1}) - (\ref{param2}) and the curvature of the domain linearly decreases with time from $\gamma = 0.4915$ to $\gamma = 0.4315$, using a value of $\varepsilon = 10^{-6}$. The initial condition used is uniformly random noise of amplitude $0.05$ for $U$ and amplitude $0.005$ for $V$. As $\gamma$ decreases further, the solution will eventually approach the black dotted curve in Figure \ref{fig:CPM_slns}.}
	\label{fig:dx_02_rand}
\end{figure}

In Figure \ref{fig:dx_02_rand} we show the results of a numerical simulation of spatial grid size $h = 0.02$ starting with an initial condition consisting of uniformly random noise with $\gamma = 0.4915$. As the system evolves and $\gamma$ slowly decreases, a $(5, 1)$ pattern emerges. This gives numerical evidence for the existence of an exponentially attracting centre manifold (\ref{EACM}).


%
%

\begin{figure}[ht]
	\centering
		\includegraphics[width=0.50\textwidth]{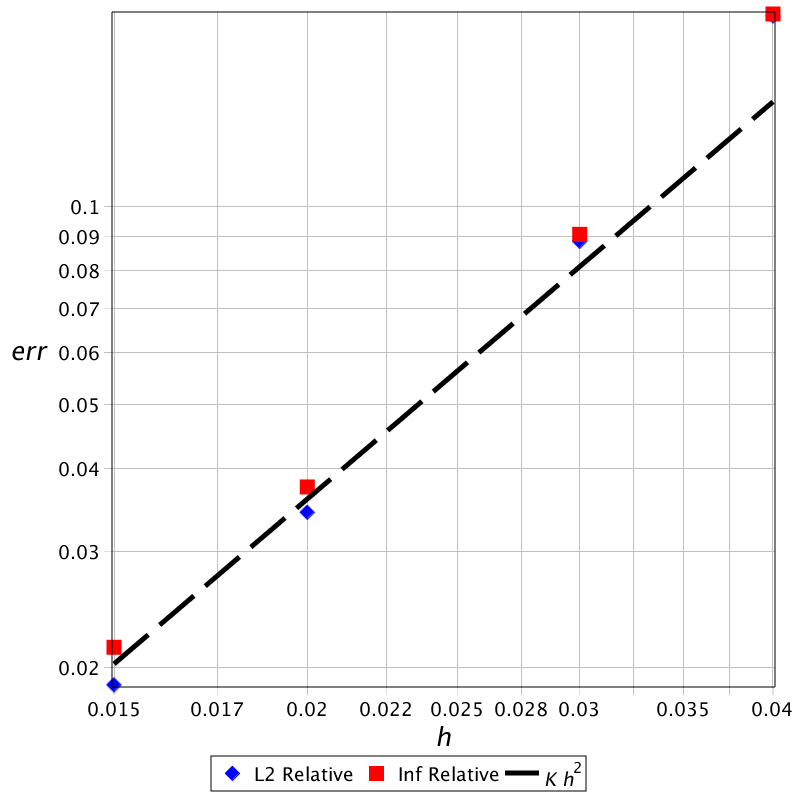}
	\caption{Convergence analysis plot for the closest point method simulations. The numerical solutions are compared with a prediction generated from the centre manifold equation. A dashed line corresponding to $h^2$ times a constant is added to compare with quadratic convergence.}
	\label{fig:Bruss_Conv3}
\end{figure}

For each simulation in Figure \ref{fig:CPM_slns} we also computed a prediction from the normal form equations using the leading term of the centre manifold (\ref{centremanifold}) as well as a truncated series of the quadratic terms for the dominant modes for the final curvature value of $\gamma = 0.4515$. We then computed an error approximation by comparing the numerical solution to this prediction for the different spatial grid sizes $h$. Figure \ref{fig:Bruss_Conv3} shows the differences compared to the grid size $h$ superimposed with a reference quadratic function to show that we have quadratic convergence in both $L^2$ and $L^{\infty}$ norms, the latter implying pointwise convergence. This is expected convergence behaviour of the closest point method \cite{RuuthMerriman2008}, so this result validates both the numerical simulations and normal form results.

\section{Conclusion}

In this work we reduced the Brusselator reaction-diffusion system (\ref{ConXvec}) on a time dependent non-isotropically slowly evolving spherical cap domain undergoing a delayed bifurcation. We did this using an asymptotic expansion of the normal form obtained from the constant (time independent) domain case. We use concepts from bifurcation theory such as centre manifold expansions and equivariant normal forms. This allows for the reaction-diffusion system to be approximated by a nonautonomous ordinary differential equation. We obtained normal forms for $m_0 \neq 0$ modes and performed calculations using the $m_0 = 5$, $n_0 = 1$ case. These reproduce the transition from the quasi-patternless state to the patterned state using significantly less computing resources compared to a numerical simulation of the partial differential equation system. This gives us a technique to better understand emergence of patterned solutions around critical stability parameter values on an evolving domain.

The existence and properties of exponentially attracting centre manifolds for systems of nonautonomous partial differential equations are well established (e.g. \cite{Chicone1997}). The centre manifold expansion we use is from \cite{Potzsche2006}, adapted to slowly changing coefficients. It remains to prove, in future work, that our formal results using the WKB method (\ref{WKB}) indeed imply the exponential dichotomy properties required for the existence of a centre manifold. Instead of a proof we provide numerical results in Section \ref{sec:Numerical} that give evidence toward its existence.

There are, however some limitations to this method. First for the asymptotic expansion to be valid we need slow domain evolution. Also, we need to start from the quasi-patternless state and go through a $m_0 \neq 0$ stability curve. More work would be required to study the $m_0 = 0$ case or transitions between patterned states. The truncation of terms in the centre manifold expansion also limits the validity of predictions to solutions close to the quasi-patternless solution. Finally, obtaining the order-$\varepsilon$  term for the cubic term $C(\varepsilon t, \varepsilon)$ of the normal form (\ref{1DNF}) could produce a better normal form approximation and should be achievable through similar methods, but with considerably more work. 

Other formulations for the curvature function may be used, for example $\gamma(\tau) = \gamma_0 - \frac{\Delta \gamma}{2}\arctan \tau$, for $-\infty < \tau < \infty,$ that has a similar, almost linear, evolution near $\tau = 0$ and asymptotically approaches constant quantities when $t$ approaches either infinities. The latter formulation could be useful in future simulations where we wish to give time for the patterns to develop in order to observe them, for example when using quicker domain changes, and for investigating exponential dichotomies, which requires at least semi-infinite time intervals.

Other future work should include more numerical simulations in order to study more transitions between different modes or the quasi-patternless solution to other modes. We could account for different ways to change the domain. Previous results in \cite{Crampin1999}, \cite{Madz2016} show that domain growth might result in more robust patterns under some circumstances. Trying to achieve this using a changing spherical cap would provide a good example with a non-isotropically changing domain.

Here we have studied a particular problem of the Brusselator system on an evolving spherical cap. However, there are morphogenesis models that use a domain whose evolution is dictated by the solution of a similar system, but incorporating biomechanical concepts. For example Brinkmann et al. \cite{Brinkmann2018} have presented a tissue surface model in three dimensions with small thickness and can show multiple morphogenesis patterns seen in many organisms. In the example of the conifer embryo, for example, we could add some thickness to the exterior tissue, or expect domain growth to be promoted in the regions of maximal concentration of one of the chemical agents. Any such features added to the model would make the asymptotic analysis much more complicated, but interesting numerical experiments could be made.


\appendix

\section{Appendix: Time derivative of the Laplace-Beltrami eigenfunction}

In the expansion of the WKB solution (\ref{WKBlin}) we are required to use the time derivative of the harmonic function $\Phi_{mn}(\tau)$. As the domain evolves with time so does the Legendre functions associated to the surfaces. This means that any point on the surface evolving in with velocity vector $\mathbf{v}$ will have a slightly different value for the function $\Phi_{mn}$. This is the time derivative we need to use.

Using domain evolution along toroidal coordinates, this corresponds with an evolving $\xi$ while keeping $\eta$ and $\phi$ constant. We can express $\Phi_{mn}$ as
\begin{equation}
\Phi_{mn}(\tau) = \cos(m (\phi + \phi_0))P^m_{\lambda_{mn}(\gamma(\tau))}(\zeta(\tau)),
\end{equation}
where $\phi_0$ is a phase parameter that we will choose, without loss of generality, to be zero. We can also express this with complex numbers
\begin{equation}
\Phi_{mn}(\tau) = e^{im (\phi + \phi_0)}P^m_{\lambda_{mn}(\gamma(\tau))}(\zeta(\tau)),
\end{equation}
which will be useful when using the symmetries of the equation. The Legendre function argument $\zeta(\tau)$ may be expressed as
\begin{align}
\zeta(\tau) &= \cos(\theta) & \text{in}& \text{ spherical coordinates,}\\
\zeta(\tau) &= \frac{1 - \cosh \eta \cos \xi(\tau)}{\cosh \eta - \cos \xi(\tau)} & \text{in}& \text{ toroidal coordinates.}
\end{align}
We can then find the time derivative using the chain rule
\begin{align}
\Phi'_{mn}(\tau) &= \cos(m (\phi + \phi_0))\left[ \frac{\partial}{\partial \lambda} P^m_{\lambda}(\zeta(\tau)) \frac{\partial \lambda_{mn}(\tau)}{\partial \tau} + \frac{\partial}{\partial \zeta} P^m_{\lambda_{mn}(\gamma(\tau))}(\zeta) \frac{\partial \zeta(\tau)}{\partial \tau} \right].
\end{align}
The derivative is thus split in two terms: the $\lambda$ derivative and the $\zeta$ derivative.
The first term does not have a practical expression that can be easily implemented in Maple, so we used finite difference approximations to find the derivative
\begin{align}
\frac{\partial}{\partial \lambda} P^m_{\lambda}(\zeta(\tau)) &\approx \frac{P^m_{\lambda + h_1}(\zeta(\tau)) - P^m_{\lambda - h_1}(\zeta(\tau))}{2h_1}\\
\frac{\partial \lambda_{mn}(\tau)}{\partial \tau} &\approx \frac{\lambda_{mn}(\tau + h_2) - \lambda_{mn}(\tau - h_2)}{2h_2},
\end{align}
for $h_1$, $h_2$ fixed, positive and very small. Note that we were already using a numerical solver to find the $\lambda$ value.

The $\zeta$ derivative term does have a useful identity
\begin{align}
\frac{\partial}{\partial \zeta}P^m_{\lambda} &= \frac{\lambda \zeta P^m_{\lambda}(\zeta) - (\lambda + m)P^m_{\lambda - 1}(\zeta)}{\zeta^2 - 1}\\
\frac{\partial \zeta}{\partial \tau} &= \frac{\sin \xi(\tau) \sinh^2 \eta}{(\cosh \eta - \cos \xi(\tau))^2} \cdot \xi' = \frac{-\gamma'}{\gamma\sqrt{1 - \gamma^2}} \cdot \frac{\gamma^2 \sinh^2 \eta}{(\cosh \eta + \sqrt{1 - \gamma^2})^2},\label{dzeta}\\
\frac{\partial \zeta}{\partial \tau} &= \frac{-\gamma'}{\gamma\sqrt{1 - \gamma^2}} \cdot \sin^2 \theta
\end{align}
where we used the identities
\begin{align}
\begin{split}
\xi' &= \frac{\gamma'}{\cos \xi}, \quad \text{ and}\\
\cos \xi &= -\sqrt{1 - \gamma^2}.
\end{split}
\end{align}
We can then use inner products of the new functions in order to find the series expression of $\Phi'_{mn}$.

\section*{Acknowledgments}
The authors thank NSERC Canada for funding and Michael J. Ward and Brian Wetton for useful discussions.

\bibliography{bibliopWG}

\begin{thebibliography}{10}

\bibitem{Ascher1995}
{\sc U.~Ascher, S.~Ruuth, and B.~Wetton}, {\em Implicit-explicit methods for
  time-dependent partial differential equations}, SIAM Journal on Numerical
  Analysis, 32 (1995), pp.~797--823.

\bibitem{Bilsborough2011}
{\sc G.~D. Bilsborough, A.~Runions, M.~Barkoulas, H.~W. Jenkins, A.~Hasson,
  C.~Galinha, P.~Laufs, A.~Hay, P.~Prusinkiewicz, and M.~Tsiantis}, {\em Model
  for the regulation of {Arabidopsis} thaliana leaf margin development},
  Proceedings of the National Academy of Sciences, 108 (2011), pp.~3424--3429.

\bibitem{Brinkmann2018}
{\sc F.~Brinkmann, M.~Mercker, T.~Richter, and A.~Marciniak-Czochra}, {\em
  Post-{Turing} tissue pattern formation: Advent of mechanochemistry}, PLOS
  Computational Biology, 14 (2018), pp.~1--21.

\bibitem{Budrene1995}
{\sc E.~O. Budrene and H.~C. Berg}, {\em Dynamics of formation of symmetrical
  patterns by chemotactic bacteria}, Nature, 376 (1995), pp.~49--53.

\bibitem{Char2018}
{\sc L.~Charette and W.~Nagata}, {\em Bifurcation of mixed mode
  reaction–diffusion patterns in spherical caps}, International Journal of
  Bifurcation and Chaos, 28 (2018), p.~1830017.

\bibitem{Chen2013}
{\sc Y.~Chen, T.~Kolokolnikov, and D.~Zhirov}, {\em Collective behaviour of
  large number of vortices in the plane}, Proc. R. Soc. A, 469 (2013),
  p.~20130085.

\bibitem{Chicone1997}
{\sc C.~Chicone and Y.~Latushkin}, {\em Center manifolds for infinite
  dimensional nonautonomous differential equations}, Journal of Differential
  Equations, 141 (1997), pp.~356 -- 399.

\bibitem{Crampin1999}
{\sc E.~J. Crampin, E.~A. Gaffney, and P.~K. Maini}, {\em Reaction and
  diffusion on growing domains: Scenarios for robust pattern formation},
  Bulletin of Mathematical Biology, 61 (1999), pp.~1093--1120.

\bibitem{Eftimie2017}
{\sc R.~Eftimie, M.~Perez, and P.-L. Buono}, {\em Pattern formation in a
  nonlocal mathematical model for the multiple roles of the tgf-β pathway in
  tumour dynamics}, Mathematical Biosciences, 289 (2017), pp.~96 -- 115.

\bibitem{Erneux1991}
{\sc T.~Erneux, E.~L. Reiss, L.~J. Holden, and M.~Georgiou}, {\em Slow passage
  through bifurcation and limit points. asymptotic theory and applications}, in
  Dynamic Bifurcations, E.~Beno{\^i}t, ed., vol.~1493 of Lecture Notes in
  Mathematics, Berlin, Heidelberg, 1991, Springer Berlin Heidelberg,
  pp.~14--28.

\bibitem{Fujita2011}
{\sc H.~Fujita, K.~Toyokura, K.~Okada, and M.~Kawaguchi}, {\em
  Reaction-diffusion pattern in shoot apical meristem of plants}, PLoS ONE, 6
  (2011), pp.~1--13.

\bibitem{vonAd2004}
{\sc L.~G. Harrison and P.~von Aderkas}, {\em Spatially quantitative control of
  the number of cotyledons in a clonal population of somatic embryos of hybrid
  larch \emph{Larix} x \emph{leptoeuropaea}}, Annals of Botany, 93 (2004),
  pp.~423--434.

\bibitem{Holl2008}
{\sc D.~M. Holloway and L.~Harrison}, {\em Pattern selection in plants:
  Coupling chemical dynamics to surface growth in three dimensions}, Annals of
  Botany, 102 (2008), pp.~361 -- 374.

\bibitem{Holl2018}
{\sc D.~M. Holloway, I.~Rozada, and J.~J.~H. Bray}, {\em Two-stage patterning
  dynamics in conifer cotyledon whorl morphogenesis}, Annals of Botany, 121
  (2018), pp.~525--534.

\bibitem{Macdonald2011}
{\sc C.~B. Macdonald, J.~Brandman, and S.~J. Ruuth}, {\em Solving eigenvalue
  problems on curved surfaces using the closest point method}, J. Comput.
  Physics, 230 (2011), pp.~7944--7956.

\bibitem{Madz2016}
{\sc A.~Madzvamuse and A.~H.~W. Chung}, {\em The bulk-surface finite element
  method for reaction–diffusion systems on stationary volumes}, Finite
  Elements in Analysis and Design, 108 (2016), pp.~9 -- 21.

\bibitem{Mein2003}
{\sc H.~Meinhardt, P.~Prusinkiewicz, and D.~Fowler}, {\em The Algorithmic
  Beauty of Sea Shells}, Virtual laboratory, Springer, 2003.

\bibitem{Nag2013}
{\sc W.~Nagata, H.~R.~Z. Zangeneh, and D.~M. Holloway}, {\em Reaction-diffusion
  patterns in plant tip morphogenesis: Bifurcations on spherical caps},
  Bulletin of Mathematical Biology, 75 (2013), pp.~2346--2371.

\bibitem{Plaza2004}
{\sc R.~Plaza, F.~S\'anchez-Gardu\~no, P.~Padilla, R.~Barrio, and P.~Maini},
  {\em The effect of growth and curvature on pattern formation}, Journal of
  Dynamics and Dif\-fer\-en\-tial Equations, 16 (2004), pp.~1093--1121.

\bibitem{Potzsche2006}
{\sc C.~P{\"o}tzsche and M.~Rasmussen}, {\em Taylor approximation of integral
  manifolds}, Journal of Dynamics and Differential Equations, 18 (2006),
  pp.~427--460.

\bibitem{Prig1968}
{\sc I.~Prigogine and R.~Lefever}, {\em Symmetry breaking instabilities in
  dissipative systems. ii}, The Journal of Chemical Physics, 48 (1968),
  pp.~1695--1700.

\bibitem{Reyn1903}
{\sc O.~Reynolds, A.~W. Brightmore, and W.~H. Moorby}, {\em Papers on
  mechanical and physical subjects}, vol.~3: The Sub-Mechanics of the Universe,
  Cambridge: The University Press, 1903.

\bibitem{RuuthMerriman2008}
{\sc S.~Ruuth and B.~Merriman}, {\em A simple embedding method for solving
  partial differential equations on surfaces}, Journal of Computational
  Physics, 227 (2008), pp.~1943--1961.

\bibitem{Stannard2008}
{\sc A.~Stannard}, {\em Dewetting-mediated pattern formation in nanoparticle
  assemblies}, Journal of Physics: Condensed Matter, 23 (2011), p.~083001.

\bibitem{Turing1952}
{\sc A.~M. Turing}, {\em The chemical basis of morphogenesis}, Philosophical
  Transactions of the Royal Society of London. Series B, Biological Sciences,
  237 (1952), pp.~37--72.

\bibitem{vanMourik2012}
{\sc S.~van Mourik, K.~Kaufmann, A.~D.~J. van Dijk, G.~C. Angenent, R.~M.~H.
  Merks, and J.~Molenaar}, {\em Simulation of organ patterning on the floral
  meristem using a polar auxin transport model}, PLoS ONE, 7 (2012), pp.~1--9.

\end{thebibliography}
{}
\bibliographystyle{siamplain}

\end{document}